\documentclass[11pt,reqno]{article}
\usepackage{amssymb, amsmath, amsthm, amsfonts, amscd, epsfig, subfig, gensymb}
\usepackage[dvipsnames]{xcolor}
\usepackage[utf8]{inputenc}
\usepackage[english]{babel}
\usepackage{bm}
\usepackage{bbm}
\usepackage{ulem}
\usepackage{graphicx}
\usepackage{geometry,graphicx,pgfplots, multirow}
\usepackage{xcolor}
\usepackage{placeins}

\usepackage[export]{adjustbox}
\usepackage[titletoc,toc,title]{appendix}
\usepackage[colorlinks=true,
            linkcolor=blue,
            citecolor=blue,
            urlcolor=blue]{hyperref}

\usepackage{algorithm}
\usepackage{algpseudocode}
\usepackage{mathtools}

\usepackage{afterpage}
\usepackage{comment}

\newtheorem{theorem}{Theorem}[section]
\newtheorem{cor}[theorem]{Corollary}
\newtheorem{lemma}[theorem]{Lemma}

\def\R{\mathbb{R}}

\newcommand{\Om}{\Omega}

\newcommand{\p}{\partial}
\newcommand{\ds}{\displaystyle}
\newcommand{\eqnref}[1]{(\ref {#1})}

\newcommand{\beq}{\begin{equation}}
\newcommand{\eeq}{\end{equation}}
\newcommand{\RN}[1]{%
  \textup{\uppercase\expandafter{\romannumeral#1}}%
}

\renewcommand{\t}[1]{\text{#1}}
\newcommand{\norm}[1]{\left\lVert#1\right\rVert}

\numberwithin{equation}{section}
\numberwithin{figure}{section}

\begin{document}

\title{
Spectral analysis and shape optimization of Steklov--Dirichlet eigenvalues on shell-type domains
}

\date{}

\author{
Dohyeong Han\thanks{Department of Mathematical Sciences, Korea Advanced Institute of Science and Technology, 291 Daehak-ro, Yuseong-gu, Daejeon 34141, Republic of Korea (dohyeonghan@kaist.ac.kr, woojoo.lee@kaist.ac.kr, mklim@kaist.ac.kr).}\and Woojoo Lee\footnotemark[1] \and Mikyoung Lim\footnotemark[1]
}

\maketitle

\begin{abstract}
We investigate the Steklov--Dirichlet eigenvalues on shell-type domains in $\mathbb{R}^n$, $n\geq 2$, where the Steklov and Dirichlet conditions are imposed on the outer and the inner boundaries, respectively. We introduce the Green function associated with the problem, which is well-defined for spectral parameters outside the spectrum and exhibits blow-up behavior at eigenvalues, and establish its spectral representation. We also examine its differentiability and derive its shape derivative under translations of the inner boundary. Using these results and a layer potential approach, we propose a computational framework for optimizing the eigenvalue with respect to the position of the inner boundary, incorporating an eigenvalue estimator for which a local error bound is established. Numerical experiments for various inner and outer boundary shapes demonstrate the effectiveness of the proposed approach in identifying configurations that maximize the first eigenvalue.

\end{abstract}



\tableofcontents

\section{Introduction}

\begin{figure}[t]
    \centering
    \includegraphics[width=0.3\textwidth]{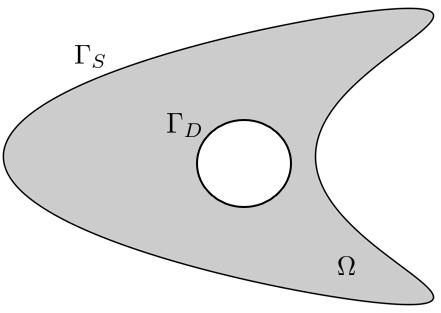}
    \caption{\label{fig:shell} A shell-type domain $\Omega$ with the boundary $\partial \Omega = \Gamma_S \cup \Gamma_D$.}
\end{figure}

We study the Steklov--Dirichlet eigenvalue problem on shell-type domains in $\mathbb{R}^n$. Specifically, let $\Omega\subset\mathbb{R}^n$ be a smooth domain whose boundary consists of two connected components, denoted by $\Gamma_D$ and $\Gamma_S$ (see Figure \ref{fig:shell}). We consider the problem
\begin{align} \label{prob:SD}
\begin{cases} 
\ds \Delta u = 0   &\ds\text{in } \Omega,\\
\ds   u = 0 & \ds\text{on }  \Gamma_D,\\
\ds   \frac{\partial u}{\partial \nu} = \sigma u &\ds\text{on } \Gamma_S
\end{cases}
\end{align}
where $\nu$ is the unit outward normal vector to $\p\Om$. A real constant $\sigma$ is said to be a Steklov--Dirichlet eigenvalue if problem \eqref{prob:SD} admits a non-trivial solution $u$. 
Assuming $\Gamma_D\neq\emptyset$, the spectrum of \eqnref{prob:SD} is discrete (see \cite{Agranovich:2006:MPS}), and the eigenvalues can be arranged in the nondecreasing sequence
$$ 0<\sigma_1(\Omega)\le \sigma_2(\Omega) \le \cdots \rightarrow \infty.$$ 

The Steklov--Dirichlet eigenvalues have been studied widely in relation to the geometry of the domain. These eigenvalues are also connected to their applications to the vibration modes of a partially free membrane, fixed along its inner boundary and with no mass in its interior \cite{Hersch:1968:EPI}, as well as to the classical Laplace eigenvalue problem in asymptotic regimes \cite{Arrieta:2008:FRL, Lamberti:2015:VSE}. As a mixed-boundary counterpart of the classical Steklov eigenvalue problem, which corresponds to the case $\Gamma_D=\emptyset$ in \eqref{prob:SD} (see \cite{Stekloff:1902:FPM} for the original formulation), the Steklov--Dirichlet problem has attracted particular attention regarding the geometric dependence of its first eigenvalue, with recent developments including \cite{Girouard:2017:SGS,Colbois:2024:SRD}. For planar doubly-connected domains, Hersch and Payne \cite{Hersch:1968:EPI} and Dittmar \cite{Dittmar:1998:IIS} established isoperimetric inequalities for the eigenvalues. Subsequently, Dittmar and Solnya obtained lower bounds for the eigenvalues \cite{Dittmar:1998:IIS, Dittmar:2005:EPC}. More recent results concern spectral stability \cite{Paoli:2021:SRS, Micetti:2022:SDS} and Riesz mean estimates \cite{Hassannezhad:2020:EBM}.

In the special case where $\Omega$ is a spherical shell, with both $\Gamma_D$ and $\Gamma_S$ being spheres, the first Steklov--Dirichlet eigenvalue can be studied in terms of the distance between the centers of the two boundary spheres. Denoting this distance by $t$, previous studies have investigated the maximality of the eigenvalue at $t=0$, which corresponds to the concentric configuration, as well as its differentiability and shape derivative with respect to $t$.
The maximality at concentric case was established in two dimensions by Seo \cite{Seo:2021:SOP} (and independently by Ftouhi \cite{Ftouhi:2022:WPS}), and in $n$ dimensions, $n\geq 3$, by Verma and Santhanam \cite{Verma:2020:EPL}. This result was also extended to convex domains with a spherical hole \cite{Gavitone:2022:IIF}. Hong, Lim, and Seo established the differentiability of the eigenvalue with respect to $t$ and derived a corresponding shape derivative \cite{Hong:2022:FSD}. Related questions concerning the shape derivatives and the dependence on $t$ have also been studied for the Dirichlet Laplacian problems \cite{Ramm:1998:IME, Chorwadwala:2013:TFL, Anoop:2018:SMF, Rane:2019:FDE} and the Zaremba problems \cite{Anoop:2021:SVR}. In addition, approximation methods for the eigenvalue on annuli ($n=2$) \cite{Hong:2022:FSD} and spherical shells in higher dimensions ($n\geq 3$) \cite{Hong:2025:AFS} have been developed. The numerical results in these works illustrate its behavior as the distance varies, such as the monotonicity pattern. 

For general shell-type domains, despite these developments, the optimization of the first Steklov--Dirichlet eigenvalue remains largely unexplored. Unlike the case of spherical shells, where the concentric configuration is known to be optimal, it is generally unclear how the inner boundary should be positioned to maximize the first eigenvalue. For instance, numerical experiments for a domain with a kite-shaped outer boundary (see Figure \ref{fig:shell}) suggest that the maximizing position lies toward the right-hand side of the domain. This example illustrates the nontrivial dependence of the optimal position on the geometric configuration of the domain and motivates the development of a shape optimization approach and a numerical algorithm for searching for optimal configurations.

In this paper, we investigate the spectral analysis of Steklov--Dirichlet eigenvalues on smooth shell-type domains and their optimization with respect to the shape of the domain. We develop a Green function-based approach that exploits the singular behavior of the Green function near eigenvalues to guide the shape optimization of the eigenvalue. This approach provides a framework for investigating how the eigenvalue depends on the position of the inner boundary and identifying configurations that maximize it. In particular, we exploit the spectral behavior of the Green function and the shape derivative to guide the optimization process with an eigenvalue estimator.


To examine the problem, motivated by \cite{Ammari:2020:OSE}, we introduce the associated Steklov--Dirichlet Green function and examine its analytic structure. We reformulate the eigenvalue problem in terms of the behavior of the Green function with respect to the spectral parameter $\lambda$. The Green function is well-defined away from the spectrum and exhibits singular behavior as the parameter approaches an eigenvalue, providing an alternative characterization of the spectrum. We establish a spectral representation of the Green function and derive its shape derivative under translations of the inner boundary, which supports the subsequent shape optimization analysis.  

Building on this spectral analysis, we develop a computational framework for determining the position of the inner boundary $\Gamma_D$ that maximizes the first eigenvalue. Rather than directly solving the eigenvalue problem, the proposed approach exploits the singular behavior of the Green function to identify eigenvalues, leading to a computationally efficient procedure. Specifically, we use the shape derivative of the Green function to guide the optimization process, avoiding the need to repeatedly compute the corresponding eigenpair. We also employ an eigenvalue estimator based on the spectral representation of the Green function, for which a local error bound is established. 
The framework is formulated using layer potential techniques and applies to general shell-type domains with the boundary $\partial\Omega=\Gamma_D\cup\Gamma_S$. 
Numerical experiments on various shell-type geometries validates the effectiveness of the proposed approach and identify candidate configurations that maximize the first eigenvalue.

The remainder of this paper is structured as follows. Section \ref{sec:prop} reviews some basic properties of the first Steklov--Dirichlet eigenvalue, arising from its variational characterization. In Section \ref{sec:spectral}, we introduce the Steklov--Dirichlet Green function, and establish its existence and derive its spectral decomposition. Subsequently, in Section \ref{sec:derv}, we discuss its differentiability and establish the shape derivative formula. Building on these results, Section \ref{sec:optm} develops a shape optimization framework associated with the eigenvalue. In Section \ref{sec:num}, we propose the computational algorithm and conduct numerical simulations for various configurations. Lastly, we conclude the paper with a brief discussion in Section \ref{sec:concl}.

\section{Variational characterization of the Steklov--Dirichlet eigenvalues} \label{sec:prop}
We recall the variational characterization of the eigenvalues\cite{Bandle:1980:IIA,Colbois:2024:SRD}: 
\begin{align}\label{variational_char}
\sigma_k=\min_{E \in \mathcal{E}(k)} \max_{0 \neq v \in E} \frac{\int_\Omega |\nabla v|^2 dx}{\int_{\Gamma_S} v^2 dS},
\end{align}
where $\mathcal{E}(k)$ consists of all $k$-dimensional subspaces of $H:=\{f\in H^1(\Omega) : f=0 \text{ on } \Gamma_D \}$.
In particular, the first eigenvalue admits the characterization
\begin{align}
    \sigma_1=\inf \left\{ \int_{\Omega} |\nabla v|^2 dx:v\in H^1(\Omega)\setminus \{0\}, v=0 \text{ on } \Gamma_D\text{, and } \int_{\Gamma_S} v^2 dS = 1 \right\}.
\end{align}
The simplicity of the first eigenvalue follows from the variational characterization as shown in \cite{Hong:2022:FSD}.

We next introduce notation in order to formulate the differentiability
result more precisely.
Let $D_1$ and $D_2$ denote the bounded regions enclosed by $\Gamma_D$ and $\Gamma_S$,
respectively, and let $d$ be a fixed unit vector in $\R^n$.
For $t \in \R$ we denote by $\Gamma_D^t$ and $D_1^t$ the translate of $\Gamma_D$ and $D_1$ by the vector $t d$.
Define $I := \{ t \in \R : \overline{D_1^t} \subset D_2 \}$ and let $I_0$ be the connected component of $I$ containing $0$.

Set $r_m := \inf I_0, \ r_M := \sup I_0$ and $\Omega^t=D_2\setminus\overline{D^t_1}$.
For each $t \in (r_m,r_M)$, we consider
\[
\begin{cases}
  \Delta u^t = 0 & \text{in } \Omega^t,\\
  u^t = 0        & \text{on } \Gamma_D^t,\\[2mm]
  \dfrac{\partial u^t}{\partial \nu} = \sigma^t u^t
                 & \text{on } \Gamma_S,
\end{cases}
\]
where $\nu$ denotes the unit outward normal vector on
$\Gamma_S$.

Choose a smooth variation field $V:\overline{D_2}\to\mathbb{R}^n$ such that
\begin{align}\label{variation_field}
V(x) = d \quad\text{for near } \overline{D^{0}_1} \quad\t{ and }
\quad
\t{supp} (V) \subset \subset D_2.    
\end{align}
We now define a map $T : (r_m,r_M)\times D_2 \to \R^n$ by
\begin{align} \label{transformation}
    T(s,x)=x+sV(x).
\end{align}

For a pair of distinct eigenvalues $(\sigma_i, \sigma_j)$, we have the following relation of the corresponding eigenfunctions,
\begin{align*}
\int_{\partial \Omega} (\partial_\nu u_i) u_j \ dS-\int_\Omega (\Delta u_i) u_j \ dx&=\int_{\partial \Omega} (\partial_\nu u_j) u_i \ dS-\int_\Omega (\Delta u_j) u_i \ dx\\
\sigma_i \int_{\Gamma_S} u_i u_j \ dS&= \sigma_j \int_{\Gamma_S} u_iu_j \ dS\\
\int_{\Gamma_S} u_i u_j \ dS &=0
\end{align*}
by Green's identity.
So we can choose $\{u_k\}_{k=1}^\infty$ to be orthonormal.
Adapting the proofs of Theorems 1 in \cite{Hong:2022:FSD}, we obtain
more general results as follows.
\begin{theorem} \label{thm:prop:diff}
    Let $D^t_1$ and $D_2$ be as defined above. Let $\sigma^t_k$ and $u^t_k$ be the simple eigenvalue of the Steklov--Dirichlet eigenvalue problem above and the corresponding eigenfunction, respectively. Then the operators
$$t \mapsto \sigma^t_k \in \R \quad \t{and} \quad t \mapsto u^t_k\circ T(t,\cdot) \in H^1(\Omega^0)$$ are $C^1$ operators with respect to $t$.
\end{theorem}

The following lemma provides an alternative variational characterization
of the Steklov--Dirichlet eigenvalues. For a subset
$A\subset \overline{D_2}$, we define
\begin{align}\label{H1_A}
H^1_A(D_2)
:=
\{f\in H^1(D_2): f=0 \text{ in } A\}. 
\end{align}

\begin{lemma}\label{Lem:var:rep}
Let $\hat{W}_{n}=\{f\in H_{\Gamma_D}^1(D_2):  \int_{\Gamma_S} fu_j \ dS=0 \t{ for } j\in \mathbb{N}, j\leq n\}$. 
Then we have
$$\sigma_{k}= \min_{w \in \hat W_{k-1} \setminus \{0\}} \frac{\int_\Omega |\nabla w|^2 dx}{\int_{\Gamma_S} w^2 dS}. $$ 
\end{lemma} 

\begin{proof}
We use proof by induction. For $k=1$, the formula clearly holds by \eqref{variational_char}.
Now assume $$\sigma_n= \min_{w \in \hat{W}_{n-1} \setminus \{0\}} \frac{\int_\Omega |\nabla w|^2 dx}{\int_{\Gamma_S} w^2 dS}.$$
Let $w \in \hat{W}_{n} \setminus \{0\}$ and $E=$ span$(u_1, \dots, u_n, w)$.
For $a \in E$, write $$a=v+cw:=a_1u_1+\cdots+ a_nu_n+cw,$$ where $a_1,\dots a_n,c$ are constants. 
Then we have
\begin{align*}
    \frac{\int_\Omega |\nabla a|^2 dx}{\int_{\Gamma_S} a^2 dS}
    = \frac{\int_\Omega |\nabla v|^2 dx+\int_\Omega |\nabla cw|^2 dx+\int_\Omega 2c\nabla v\cdot\nabla w  dx}{\int_{\Gamma_S} v^2 dS+\int_{\Gamma_S} c^2w^2 dS+\int_{\Gamma_S} 2c(vw) dS}.
\end{align*}
Note that $\int_{\Gamma_S} vw \ dS=0$ by definition of $w$ and 
$$\int_\Omega \nabla v\cdot\nabla w  dx=\sum_{k=1}^n a_k\int_{\Omega} \nabla u_k \cdot \nabla w \ dx=\sum_{k=1}^n \sigma_k a_k\int_{\Gamma_S} u_k w \ dS=0.$$
Thus, we have
\begin{align*}
    \frac{\int_\Omega |\nabla a|^2 dx}{\int_{\Gamma_S} a^2 dS}
    = \frac{\int_\Omega |\nabla v|^2 dx+\int_\Omega |\nabla cw|^2 dx}{\int_{\Gamma_S} v^2 dS+\int_{\Gamma_S} c^2w^2 dS}\leq \max\left\{\frac{\int_\Omega |\nabla v|^2 dx}{\int_{\Gamma_S} v^2 dS},\frac{\int_\Omega |\nabla w|^2 dx}{\int_{\Gamma_S} w^2 dS}  \right\}.
\end{align*}
Applying the same argument to $v$, we obtain
$$\frac{\int_\Omega |\nabla a|^2 dx}{\int_{\Gamma_S} a^2 dS}
    \leq \max\left\{\sigma_n,\frac{\int_\Omega |\nabla w|^2 dx}{\int_{\Gamma_S} w^2 dS}  \right\}.$$
Since $w\in \hat{W}_{n} \setminus \{0\}$, we have
$$\sigma_n \leq\frac{\int_\Omega |\nabla w|^2 dx}{\int_{\Gamma_S} w^2 dS}$$
by the induction hypothesis.
Hence,
$$\sigma_{n+1}\leq \max_{a\in E \setminus \{0\}}\frac{\int_\Omega |\nabla a|^2 dx}{\int_{\Gamma_S} a^2 dS}
= \frac{\int_\Omega |\nabla w|^2 dx}{\int_{\Gamma_S} w^2 dS} \quad \text{for all } w\in \hat{W}_n \setminus \{ 0 \}.$$
On the other hand, for $u_{n+1}\in \hat{W}_n \setminus \{ 0 \}$,
\begin{align*}
    \frac{\int_\Omega |\nabla u_{n+1}|^2 dx}{\int_{\Gamma_S} u_{n+1}^2 dS}=\sigma_{n+1}.
\end{align*}
Therefore, we have
$$\sigma_{n+1}= \min_{w \in \hat{W}_{n} \setminus \{0\}} \frac{\int_\Omega |\nabla w|^2 dx}{\int_{\Gamma_S} w^2 dS}. $$
\end{proof}


\section{Spectral decomposition of the Green function for the Steklov--Dirichlet problem} \label{sec:spectral}
\subsection{Green functions on shell-type domains}
We introduce, for each $x_s\in\Omega$, the Green function $G_\lambda(x_s,\cdot)$ associated with the Steklov--Dirichlet boundary value problem: 
\begin{align}\label{greeneq}
    \begin{cases}
  \Delta G_\lambda(x_s,y) = -\delta_{x_s}(y) & \text{in } \Omega,\\
  G_\lambda(x_s,y) = 0        & \text{on }  \Gamma_D,\\[2mm]
  \dfrac{\partial G_\lambda}{\partial \nu}(x_s,y) = \lambda G_\lambda(x_s,y)
                 & \text{on }  \Gamma_S,
\end{cases}
\end{align}
where $\lambda$ is assumed not to belong to the spectrum of the corresponding
Steklov--Dirichlet problem. Refer to the construction in \cite{Ammari:2020:OSE}.

\begin{lemma}\label{Lem:Green:exist}
    Let $g\in L^2(\Gamma_S)$.
    Then the equation
    \begin{align*}
        \begin{cases}
            \Delta U(x) = 0 & \text{in } \Omega,\\
            U(x) = 0        & \text{on } \Gamma_D,\\[2mm]
            \dfrac{\partial U}{\partial \nu}(x) = \lambda U(x) + g(x) & \text{on } \Gamma_S,
        \end{cases}
    \end{align*}
has a unique solution.
\end{lemma}
\begin{proof}
We can identify $H^1_{\Gamma_D}(\Omega)$ with $H^1_{\bar{D_1}}(D_2)$.
We have the weak formulation:
\begin{align}\label{weakf}
  \int_{D_2} \nabla U\cdot \,\nabla\varphi\,dx
  = \int_{\Gamma_S} \left(\lambda U+g\right)\varphi\,dS 
\end{align}
for all $\varphi\in H^1_{\bar{D_1}}(D_2)$.
We define two compact linear maps (see the proof of theorem 1 in \cite{Hong:2022:FSD}) $S_1, S_2 : H^{1}_{\bar{D}_1}(D_2) \rightarrow \left(H^{1}_{\bar{D}_1}(D_2)\right)'$ by
$$
\langle S_1(v), \psi \rangle = \int_{\Gamma_S} v\psi \, dS,
$$
$$
\langle S_2(v), \psi \rangle = \int_{D_2} v\psi \, dx \qquad \text{for } v, \psi \in H^{1}_{\bar{D}_1}(D_2).
$$
We can rewrite \eqref{weakf} as
\begin{align}
    \int_{D_2} \nabla U\cdot \,\nabla\varphi\,dx
  - \int_{\Gamma_S} \left(\lambda U\right)\varphi\,dS=(U,\varphi)- \left\langle (S_2+\lambda S_1)(U), \varphi \right\rangle =\int_{\Gamma_S} g \varphi \ dS
\end{align}
where $(\cdot,\cdot)$ is the inner product on $H^1_{\bar{D_1}}(D_2)$.
Note that $g\in L^2(\Gamma_S)$ and $$\langle g, \varphi \rangle =\int_{\Gamma_S} g \varphi dS \in \left(H^{1}_{\bar{D}_1}(D_2)\right)'$$
by trace theorem.
Then by Fredholm alternative, we have a unique solution $U$ up to $\ker{(Id-S_2-\lambda S_1)}$ if and only if 
$\langle g, \cdot \rangle \perp \ker{(Id-S_2-\lambda S_1)}$. Since $\lambda$ is not an eigenvalue of corresponding Steklov--Dirichlet problem, we have a trivial kernel. Hence, the solution $U$ exists uniquely.
\end{proof}
Using Lemma~\ref{Lem:var:rep} and Lemma~\ref{Lem:Green:exist}, we establish the existence of the Green function and derive its spectral decomposition in the following subsection.

\subsection{Spectral decomposition of the Green function}

In this subsection, we derive spectral representations of the Green function in the interior of the domain and of its trace on the outer boundary.

\begin{theorem} \label{Greendef:thm:diff}
For all $x_s\in \Omega$ and $\lambda$ that is not in the spectrum of the problem, $G_\lambda(x_s,\cdot)$ exists and $G_\lambda(x_s,\cdot)\in L^p(\Omega)$ where $1 \leq p < \frac{n}{n-2}$ for $n\geq 3$ and $1 \leq p < \infty$ for $n=2$. Furthermore, 
$$G_\lambda(x_s,y)=\Phi(x_s,y) + \sum_{k=1}^\infty \frac{u_k(x_s)u_k(y)}{\sigma_k-\lambda}+\mathcal{R}_{x_s}(y),$$
where $\mathcal{R}_{x_s} \in H^1(\Omega) \cap C^\infty(\Omega).$
\end{theorem}
\begin{proof}
We use the decomposition
$$G_\lambda(y)=\Phi(x_s,y)+Z(y)+U(y),$$ 
where $\Phi$ is the fundamental solution of the Laplacian, and $Z$ and $U$ solve
\begin{align*}
    \begin{cases}
        \Delta Z=0 \quad \t{in }\Omega,\\
        Z=0  \quad \t{on } \Gamma_S,\\
        Z=-\Phi  \quad \t{on } \Gamma_D,
    \end{cases}
    \qquad \qquad
    \begin{cases}
        \Delta U=0 \quad \t{in} \Omega,\\
        \dfrac{\partial U}{\partial \nu}=\lambda U-\partial_\nu \Phi +\lambda \Phi -\partial_\nu Z  \quad \t{on } \Gamma_S,\\
        U=0  \quad \t{on } \Gamma_D.
    \end{cases}
\end{align*}
Since $Z$ is a solution of the Dirichlet problem, such $Z$ uniquely exists. $U$ also uniquely exists by Lemma \ref{Lem:Green:exist}, since $-\partial_\nu \Phi +\lambda \Phi -\partial_\nu Z\in L^2(\Gamma_S)$. Thus, $G_\lambda$ uniquely exists with $G_\lambda\in L^p(\Omega)$. 

Let $W=\{v \in H^1(\Omega): v \t{ is a weak solution of } \Delta v=0 \t{ and } v=0 \t{ on } \Gamma_D\}$. For $j\geq 1$, define $W_j=$ span$(u_1, \dots, u_j)$ and $K_j=W_j^\perp$. Here we consider an inner product $(f,g)_W=\int_{\Omega} \nabla f \cdot \nabla g \, dx$ and an orthonormal set $A=\{\frac{u_k}{\sqrt{\sigma_k}} \}_{k=1}^\infty$.
Suppose $f\in \cap_{j=1}^\infty K_j\setminus \{0\}$. Then by lemma \ref{Lem:var:rep},
\begin{align*}
    \sigma_j \leq \frac{\int_\Omega |\nabla f|^2 dx}{\int_{\Gamma_S} f^2 dS}
\end{align*}
for all $j\in \mathbb{N}$. Since $\sigma_j \to \infty$, we have a contradiction. Thus, $A$ is an orthonormal basis of $W$.
So we can write
\begin{align*}
    U(y)=\sum_{k=1}^\infty (U,e_k)_W \ e_k(y),
\end{align*}
where $e_k = \frac{u_k}{\sqrt{\sigma_k}}$.
Let $(f,g)_{\Gamma_S}=\int_{\Gamma_S} fg \, dS$. 
Then we have $(U,e_i)_W=\sqrt{\sigma_i} \, (U,u_i)_{\Gamma_S}$ and
$$ U(y)=\sum_{k=1}^\infty (U,u_k)_{\Gamma_S} \ u_k(y).$$

Using Green's identity, we have
\begin{align} \label{G_co}
    \int_{\Gamma_S}G_\lambda u_k \, dS=\frac{u_k(x_s)}{\sigma_k-\lambda}.
\end{align}
Note that 
$$\int_{\Gamma_S}G_\lambda u_k \, dS=\int_{\Gamma_S}(U+\Phi) u_k \, dS. $$
Therefore, we have the coefficient,
$$(\Phi,u_k)_{\Gamma_S}=\frac{u_k(x_s)}{\sigma_k-\lambda}-(U,u_k)_{\Gamma_S}.$$
Now define $P_{x_s}$ as a solution of the equation
\begin{align*}
    \begin{cases}
        \Delta P_{x_s} = 0 \quad \text{in } \Omega, \\
        P_{x_s} = 0 \quad \text{on } \Gamma_D,\\
        P_{x_s} = \Phi(x_s,\cdot) \quad \text{on } \Gamma_S.
    \end{cases}
\end{align*}
Then we have 
$$P_{x_s}=\sum^\infty_{k=1} (P_{x_s},e_k)_W \, e_k = \sum^\infty_{k=1} (\Phi,u_k)_{\Gamma_S}u_k \in H^1(\Omega)\cap C^\infty(\Omega).$$
Therefore,
$$\mathcal{R}_{x_s}(y):= Z(y)-P_{x_s}(y) = Z(y)+U(y)-\sum_{k=1}^\infty \frac{u_k(x_s)u_k(y)}{\sigma_k-\lambda}\in H^1(\Omega)\cap C^\infty(\Omega),$$
and 
$$G_\lambda(x_s,y)=\Phi(x_s,y) + \sum_{k=1}^\infty \frac{u_k(x_s)u_k(y)}{\sigma_k-\lambda}+\mathcal{R}_{x_s}(y).$$
\end{proof}

\begin{theorem}\label{thm:Numerical:J_expression}
For every $\lambda$ outside the Steklov--Dirichlet spectrum, we have
\[
    \int_{\Gamma_S}|G_\lambda(x_s,x)|^2\,dS_x
    =
    \sum_{j=1}^{\infty}
    \frac{|u_j(x_s)|^2}{(\sigma_j-\lambda)^2}.
\]

\end{theorem}

\begin{proof}
The traces $\{u_k|_{\Gamma_S}\}_{k=1}^{\infty}$ form an
orthonormal basis of $L^2(\Gamma_S)$. Orthonormality follows from
our normalization. To prove completeness, let
$f\in H^{1/2}(\Gamma_S)$ and let $v_f\in W$ be the unique harmonic
extension satisfying
\[
    v_f=0 \quad\text{on }\Gamma_D,
    \qquad
    v_f=f \quad\text{on }\Gamma_S.
\]
Such an extension exists by the well-posedness of the Dirichlet
problem. Green's identity gives
\[
    (v_f,e_k)_W
    =\sqrt{\sigma_k}\,(f,u_k)_{\Gamma_S}.
\]
Since $\{e_k\}_{k=1}^{\infty}$ is an orthonormal basis of $W$
and the energy norm on $W$ is equivalent to the $H^1(\Omega)$-norm,
we obtain
\[
    v_f
    =\sum_{k=1}^{\infty}(f,u_k)_{\Gamma_S}\,u_k
    \qquad\text{in }H^1(\Omega).
\]
Taking traces yields
\[
    f
    =\sum_{k=1}^{\infty}(f,u_k)_{\Gamma_S}\,
      u_k|_{\Gamma_S}
    \qquad\text{in }L^2(\Gamma_S).
\]
Thus, the $L^2(\Gamma_S)$-closure of the span of these traces
contains $H^{1/2}(\Gamma_S)$. Since $H^{1/2}(\Gamma_S)$ is dense
in $L^2(\Gamma_S)$, completeness follows.
    Since $\{u_j|_{\Gamma_S}\}_{j=1}^{\infty}$ is an
orthogonal basis of $L^2(\Gamma_S)$, we may expand the trace
of the Green function as
\[
    G_\lambda(x_s,\cdot)\big|_{\Gamma_S}
    =
    \sum_{j=1}^{\infty}
    c_j u_j\big|_{\Gamma_S}
    \qquad\text{in }L^2(\Gamma_S),
\]
where
\[
    c_j
    =
    \int_{\Gamma_S}
    G_\lambda(x_s,x)u_j(x)\,dS_x.
\]

We already have the coefficients $c_j$ in \eqref{G_co},
\[
    G_\lambda(x_s,\cdot)\big|_{\Gamma_S}
    =
    \sum_{j=1}^{\infty}
    \frac{u_j(x_s)}
    {\sigma_j-\lambda}
    u_j^q\big|_{\Gamma_S}
    \qquad\text{in }L^2(\Gamma_S).
\]
Finally, Parseval's identity gives
\begin{align*}
    \int_{\Gamma_S}|G_\lambda(x_s,x)|^2\,dS_x
    &=
    \left\|
        G_\lambda^q(x_s,\cdot)
    \right\|_{L^2(\Gamma_S)}^2 \\
    &=
    \sum_{j=1}^{\infty}
    \frac{|u_j(x_s)|^2}
    {(\sigma_j-\lambda)^2}.
\end{align*}
\end{proof}

The spectral decomposition of the Green function reveals its singular behavior as the spectral parameter approaches a Steklov--Dirichlet eigenvalue. This property plays a key role in developing our eigenvalue optimization algorithm.

\section{Shape derivative of the Green function} \label{sec:derv}
We continue our analysis of the Green function by establishing its differentiability and deriving its shape derivative. We consider the Green function $G_\lambda^t$ defined on the domain $\Omega^t = D_2 \backslash \overline{D_1^t}$ (see Section \ref{sec:prop}).

\subsection{Differentiability}
To remove the singularity of the Green function, we set 
\begin{align}\label{wdef}
    w^t(x)=G^t_\lambda(x_s,x)-\Phi(x_s,x),
\end{align}
Then $w^t$ satisfies the following boundary value problem:
\begin{align}
\begin{cases}\label{weq}
  \Delta w^t(x) = 0 & \text{in } \Omega^t,\\
  w^t(x) = -\Phi(x_s,x)        & \text{on } \Gamma_D^t,\\[2mm]
  \dfrac{\partial w^t}{\partial \nu}(x) = \lambda w^t(x)+\lambda \Phi(x_s,x)-\dfrac{\partial\Phi}{\partial \nu}(x)
                 & \text{on } \Gamma_S.
\end{cases}    
\end{align}
The solution of \eqref{weq} satisfies \eqref{wdef} and is unique, since the solution of \eqref{greeneq} is unique.

\begin{theorem}\label{thm:Greendef:diff}
    The mapping $t \mapsto w^{t}\circ T(t,\cdot) \in H^1(\Omega^0)$ is $C^1$ in a sufficiently small neighborhood of $t=0$.
\end{theorem}

\begin{proof}
We adapt the proof strategy of \cite{Hong:2022:FSD} to our setting. 
Recall that the variation field $V$ and the associated transformation $T:(r_m,r_M)\times D_2\to\mathbb R^n$, introduced in \eqref{variation_field} and \eqref{transformation}, are given by
$$
T(s,x)=x+sV(x).
$$
Note that $T(s,D_2)=D_2$ and $w^{s} \circ T(s,\cdot)\in H^1(\Omega^0)$ for sufficiently small $s$.
Since 
\begin{align}
    DT(0,\cdot)=Id, \label{eq0}
\end{align}
there is a neighborhood $U$ of $0$ such that $T(s,\cdot)$ is a diffeomorphism of $D_2$ onto itself.

We can consider any $f\in H^1(\Omega^0)$ as an element of $H^1(D_2)$ using a bounded linear extension operator (see the theorem 7.25 in \cite{Gilbarg:2001:EPS}). We fix a bounded linear extension operator $E^0:H^1(\Omega^0) \to H^1(D_2)$.
Then the following weak formulation holds:
\begin{align}
  \int_{\Omega^s} \nabla (w^{s})\cdot \,\nabla\varphi\,dx
  = \int_{\Gamma_S} \left(\lambda w^{s}+\lambda \Phi-\partial_\nu \Phi\right)\varphi\,dS \label{weak}
\end{align}
for all $\varphi\in H^1_{\overline{D_1}^{s}}(D_2)$ (see \eqref{H1_A}).
By the change of variables $y=T(s,x)$ in \eqref{weak} and using the chain rule,
we obtain
\begin{align}
\int_{D_2}
  \nabla\bigl(E^{0}(w^{s}\circ T(s,\cdot))\bigr) \bigl(DT(s,\cdot)\bigr)^{-1}
  \cdot
  \nabla\bigl(\varphi\circ T(s,\cdot)\bigr)\,\bigl(DT(s,\cdot)\bigr)^{-1}
  \,|DT(s,\cdot)|\,dx \notag \\
= \int_{\Gamma_S} \left( \lambda w^{s}+
\lambda\Phi
-\partial_\nu \Phi \right)
\varphi \, dS, \label{eq}
\end{align}
since $\varphi=0$ in $\overline{D}_1^{s}$.
\vspace{0.5em}

We denote by $(H^1_{\bar{D}_1^{0}}(D_2))'$ the dual space of $H^1_{\bar{D}_1^{0}}(D_2)$
and denote by $\langle\cdot,\cdot\rangle$ the dual pairing between $(H^1_{\bar{D}_1^{0}}(D_2))'$ and $H^1_{\bar{D}_1^{0}}(D_2)$.
We then define 
$$f=(f_1,f_2):U\times H^1(\Omega^0)\to (H^1_{\bar{D}_1^{0}}(D_2))' \times H^{1/2}(\Gamma_D)$$
by
\begin{align}
\notag
\langle f_1(s,v),\psi\rangle
=
\int_{D_2}
  \bigl((\nabla E^{0}v)(DT(s,\cdot))^{-1}\bigr)\cdot
  \bigl((\nabla\psi)(DT(s,\cdot))^{-1}\bigr)
  |DT(s,\cdot)|\,dx\\
 -\int_{\Gamma_S} \left( \lambda v+
\lambda\Phi
-\partial_\nu \Phi\right)
\psi \, dS
\end{align}
for all $\psi\in H^1_{\bar{D}_1^{0}}(D_2)$ and
\begin{align}
    f_2(s,v)=v|_{\Gamma_D}+\Phi\circ T(s,\cdot)|_{\Gamma_D}.
\end{align}
  
The mapping $f$ is $C^1$ in a neighborhood of $(0,w^{0})$.
Define
$$
g:U\to H^1(\Omega^0),\qquad
g(s) := w^{s}\circ T(s,\cdot).
$$
Then we have by \eqref{eq} and definition of $w^s$,
\[
f\bigl(s,g(s)\bigr) = 0\qquad\text{for all }s\in U.
\]

To apply the implicit function theorem, we need to show that the partial derivative
\begin{align}
\frac{\partial f}{\partial v}\bigg|_{(0,w^{0})}: H^1(\Omega^0)\to (H^1_{\bar{D}_1^{0}}(D_2))' \times H^{1/2}(\Gamma_D) \label{defpartial}
\end{align}
is an isomorphism.
That is, for any $(h,k)\in \bigl(H^1_{\overline{D}_1^{0}}(D_2)\bigr)'\times H^{1/2}(\Gamma_D)$,
we can find a unique element $\omega \in H^1(\Omega^0)$ such that
\begin{equation}
\langle h,\psi\rangle
=
\left\langle
\frac{\partial f_1}{\partial v}
\bigg|_{(0,w^{0})}
(\omega),\psi
\right\rangle,
\quad
\t{ for all }\,\psi\in H^1_{\overline{D}_1^{0}}(D_2)
\label{partialdiff}
\end{equation}
and
\begin{align}
    k=\frac{\partial f_2}{\partial v}\bigg|_{(0,w^{0})}(\omega)
\end{align}
From \eqref{eq0} and the definition of $f$, the right–hand sides are
\begin{align}
\left\langle
\frac{\partial f_1}{\partial v}
\bigg|_{(0,w^{0})}(\omega),
\psi
\right\rangle
&= 
\int_{D_2}\nabla (E^{0}\omega)\cdot\nabla\psi\,dx
-\lambda \int_{\Gamma_S}\omega\psi\,dS 
\end{align}
and
\begin{align}
    \frac{\partial f_2}{\partial v}\bigg|_{(0,w^{0})}(\omega)=\omega|_{\Gamma_D}.
\end{align}
Since $k\in H^{1/2}(\Gamma_D)$, we have an extension $\bar{k}\in H^1(\Omega^0)$. So we rewrite $\omega=\omega_0+\bar{k}$, where $\omega_0 \in H^1_{\overline{D}_1^{0}}(D_2)$. Note that we identify $H^1_{\overline{D}_1^{0}}(D_2)$ with $\{ f \in H^1(\Omega): f|_{\partial D_1^{0}}=0 \}$.\\
Then \eqref{partialdiff} becomes
\begin{align} \label{iso_condition_w0}
    \langle \tilde{h},\psi\rangle
=
\left\langle
\frac{\partial f_1}{\partial v}
\bigg|_{(0,w^{0})}
(\omega_0),\psi
\right\rangle,
\quad
\t{ for all }\,\psi\in H^1_{\overline{D}_1^{0}}(D_2),
\end{align}
$$\t{where } \langle \tilde{h},\psi\rangle=\left\langle h-\frac{\partial f_1}{\partial v}
\bigg|_{(0,w^{0})}(\bar{k}),\psi\right\rangle.$$

From now on, we consider \eqref{iso_condition_w0} and find a unique solution $\omega_0$. 
We define two compact linear maps $S_1, S_2 : H^{1}_{\bar{D}_1^{0}}(D_2) \rightarrow \left(H^{1}_{\bar{D}_1^{0}}(D_2)\right)'$ by
$$
\langle S_1(v), \psi \rangle = \int_{\Gamma_S} v\psi \, dS,
$$
$$
\langle S_2(v), \psi \rangle = \int_{D_2} v\psi \, dx \qquad \text{for } v, \psi \in H^{1}_{\bar{D}_1^{0}}(D_2).
$$

For $\omega_0\in H^1_{\bar{D}_1^{0}}(D_2)$, we can rewrite \eqref{iso_condition_w0} as
\begin{align} \label{fredholm apply}
(\omega_0, \psi) - \langle (\lambda S_1 + S_2)(\omega_0), \psi \rangle = \langle \tilde{h},\psi \rangle    
\end{align}
where $(\omega_0, \psi)$ is the inner product in $H^{1}_{\bar{D}_1^{0}}(D_2)$.
We can regard $\tilde{h}$ and $(\lambda S_1 + S_2)(\omega_0)$ as elements in the dual space of
$H^{1}_{\bar{D}_1^{0}}(D_2)$. Since $\lambda S_1 + S_2$ is compact and
self-adjoint, there is a unique solution $\omega_0\in H^1_{\bar{D}_1^{0}}(D_2)$ to \eqref{fredholm apply} up to 
$$\ker(\mathrm{Id} - (\lambda S_1 + S_2)),$$ 
if and only if  
$$\tilde{h} \perp \ker(\mathrm{Id} - (\lambda S_1 + S_2))$$
by the Fredholm alternative. 
Note that $\ker(\mathrm{Id} - (\lambda S_1 + S_2))$ is the eigenspace of $\lambda$ of the original Steklov--Dirichlet problem and $\lambda$ is not the eigenvalue. Thus we obtain the trivial kernel and a unique solution $\omega_0\in
H^1_{\overline{D}_1^{0}}(D_2)$. So we have a unique solution $\omega=\omega_0+\bar{k}\in H^1(\Omega^0)$.
Therefore, the partial derivative \eqref{defpartial} is an isomorphism and this finishes the proof by the implicit function theorem.
\end{proof}

We consider $w^t\circ T(t,\cdot)$ as an element in $H^1(D_2)$, using extension
\begin{align*}
    \bar{w}^t(x)=\tilde{w}^t\circ T(t,x), \t{ where } \tilde{w}^t=
    \begin{cases}
        w^t(x), \quad \t{ for }x\in \Omega^t, \\
        -\Phi(x), \quad \t{ for }x\in \overline{D}_1^t.
    \end{cases}
\end{align*}

\begin{cor} \label{Cor:green:H1diff}
    The mapping $t \mapsto \bar{w}^t \in H^1(D_2)$ is differentiable with respect to $t$.
\end{cor}
    
\begin{proof}
We need to find a function $F\in H^1(D_2)$ such that
\begin{align}
    \lim_{h\to 0} \left\Vert\frac{(\bar{w}^{t_0+h}-\bar{w}^{t_0})(\cdot)}{h}-F(\cdot)\right\Vert_{H^1(D_2)}=0.
\end{align}
Observe
\begin{align*}
    &\left\Vert\frac{(\bar{w}^{t_0+h}-\bar{w}^{t_0})(\cdot)}{h}-F(\cdot)\right\Vert_{{H^1}(D_2)}\\
    \leq &\left\Vert\frac{(\bar{w}^{t_0+h}-\bar{w}^{t_0})(\cdot)}{h}-F(\cdot)\right\Vert_{H^1(\Omega^0)}+
    \left\Vert\frac{(\bar{w}^{t_0+h}-\bar{w}^{t_0})(\cdot)}{h}-F(\cdot)\right\Vert_{H^1(D_1^0)}.
\end{align*}
By the differentiability on $H^1(\Omega^0)$ (see Theorem \ref{thm:Greendef:diff}), we have $F_1\in H^1(\Omega^0)$ such that
\begin{align}
    \lim_{h\to 0}\left\Vert\frac{(\bar{w}^{t_0+h}-\bar{w}^{t_0})(\cdot)}{h}-F_1(\cdot)\right\Vert_{H^1(\Omega^0)}=0.
\end{align}
By the trace theorem, we also have
\begin{align*}
    \left\Vert\frac{(\bar{w}^{t_0+h}-\bar{w}^{t_0})(\cdot)}{h}-F_1(\cdot)\right\Vert_{L^2(\Gamma_D^0)}\leq C\left\Vert\frac{(\bar{w}^{t_0+h}-\bar{w}^{t_0})(\cdot)}{h}-F_1(\cdot)\right\Vert_{H^1(\Omega^0)}.
\end{align*}
Therefore,
\begin{align*}
&\lim_{h\to0} \left\Vert\frac{(\bar{w}^{t_0+h}-\bar{w}^{t_0})(\cdot)}{h}-F_1(\cdot)\right\Vert_{L^2(\Gamma_D)}=0.
\end{align*}
We know the exact trace of $\bar{w}^{t_0+h}$ and the difference quotients converge in $L^2$.  
Thus, the trace of $F_1$ on $\Gamma_D$ is $-\nabla\Phi(\cdot+t_0d)\cdot d$.

On the other hand, we clearly have
\begin{align*}
\lim_{h\to0}\left\Vert\frac{(\bar{w}^{t_0+h}-\bar{w}^{t_0})(\cdot)}{h}+\nabla\Phi(\cdot+t_0d)\cdot d\right\Vert_{H^1(D_1^0)}=0.
\end{align*}
Lastly, we define $F\in H^1(D_2)$ as follows:
\begin{align*}
    F(x)=
    \begin{cases}
        F_1(x), \quad \t{ for }x\in \Omega^0, \\
        -\nabla\Phi(x+t_0d)\cdot d, \quad \t{ for }x\in \overline{D^0_1}.
    \end{cases}
\end{align*}
\end{proof}

\begin{cor} \label{Cor:Green:l2diff}
    The mapping $t \mapsto \tilde{w}^t\in L^2(D_2)$ is differentiable at $t=0$. In particular, the mapping $t \mapsto G^t_\lambda-\Phi$ is also differentiable at $t=0$.
\end{cor}

\begin{proof}
We follow the decomposition used in the proof of
\cite[Lemma~5.2.7]{Henrot:2005:VOF}, working in $L^2$.
Let
\[
    \mathcal E:H^1(D_2)\longrightarrow H^1(\mathbb R^n)
\]
be a bounded linear extension operator, and set
\[
    h_t:=\mathcal E\bar w^t.
\]
Extend $V\in C_c^\infty(D_2;\mathbb R^n)$ by zero to
$\mathbb R^n$, retaining the same notation, and write
\[
    T_t:=\operatorname{Id}+tV,
    \qquad
    \psi_t:=T_t^{-1}.
\]
For sufficiently small $|t|$, $T_t$ is a diffeomorphism
of $\mathbb R^n$ preserving $D_2$.
By Corollary~\ref{Cor:green:H1diff}, there exists
$\dot h_0\in H^1(\mathbb R^n)$ such that
\[
    \frac{h_t-h_0}{t}\longrightarrow\dot h_0
    \qquad\text{in }H^1(\mathbb R^n).
\]

For each fixed $f\in H^1(\mathbb R^n)$, define
\[
    \mathcal C_f(\theta)
    :=f\circ(\operatorname{Id}+\theta)^{-1}
\]
for $\theta$ in a sufficiently small neighborhood of $0$ in
$W^{1,\infty}(\mathbb R^n;\mathbb R^n)$.
The generalized part of
\cite[Lemma~5.3.3]{Henrot:2005:VOF}, applied with $p=2$,
$g(\theta)\equiv f$, and
$\Psi(\theta)=(\operatorname{Id}+\theta)^{-1}$,
together with the inverse-map estimates in
\cite[(5.5)]{Henrot:2005:VOF},
shows that $\mathcal C_f$ is differentiable at $0$
as an $L^2(\mathbb R^n)$-valued mapping.
Since $t\mapsto tV$ is bounded and linear, it follows that
\[
    t\longmapsto\mathcal C_f(tV)=f\circ\psi_t
\]
is differentiable at $t=0$ in $L^2(\mathbb R^n)$.

We now write
\[
\begin{aligned}
    \frac{h_t\circ\psi_t-h_0}{t}
    &=
    \left(\frac{h_t-h_0}{t}-\dot h_0\right)\circ\psi_t
    +\dot h_0\circ\psi_t\\
    &\quad+
    \frac{h_0\circ\psi_t-h_0}{t}.
\end{aligned}
\]
The first term tends to zero in $L^2(\mathbb R^n)$ by the
uniform change-of-variables estimate
\[
    \|F\circ\psi_t\|_{L^2(\mathbb R^n)}
    \le C\|F\|_{L^2(\mathbb R^n)},
\]
where $C$ is independent of sufficiently small $t$.
The second term tends to $\dot h_0$ by the continuity of
$\mathcal C_{\dot h_0}$ at $0$.
The third term converges by the differentiability of
$t\mapsto\mathcal C_{h_0}(tV)$ at $0$.
Consequently, $t\mapsto h_t\circ\psi_t$ is differentiable
at $t=0$ in $L^2(\mathbb R^n)$.

Finally, since $\psi_t(D_2)=D_2$, we have
\[
    \left.(h_t\circ\psi_t)\right|_{D_2}
    =\bar w^t\circ T_t^{-1}
    =\widetilde w^t.
\]
Restriction to $D_2$ therefore proves the differentiability
of $t\mapsto\widetilde w^t$ at $t=0$ in $L^2(D_2)$.
\end{proof}

Theorem \ref{thm:Greendef:diff} establishes the $H^1$-differentiability of $w^t\circ T(t,\cdot)$ on $\Omega^0$, while Corollary \ref{Cor:green:H1diff} proves the $H^1$-differentiability of $\bar{w}^t$ on $D_2$ and Corollary \ref{Cor:Green:l2diff} derives the $L^2$-differentiability of $\tilde{w}^t$ on $D_2$. These differentiability results then allow us to obtain the governing equation for shape derivative of $G_\lambda^t$ at $t=0$.

\subsection{Shape derivative formula}
We denote shape derivative of $G_\lambda^t$ at $t=0$ by $G_\lambda'$.
\begin{theorem} \label{thm:derv:G'} 
The derivative of the mapping $t \mapsto G_\lambda^t$ at $t=0$ satisfies the following equation,
\begin{align*}
    \begin{cases}
        \Delta G_\lambda'(x) = 0 & \text{in } \Omega^{0},\\
  G_\lambda'(x) = - \frac{\partial G^0_\lambda}{\partial \nu} (\nu \cdot d)       & \text{on } \Gamma_D^{0},\\[2mm]
  \dfrac{\partial G_\lambda'}{\partial \nu}(x) = \lambda G_\lambda'
                 & \text{on } \Gamma_S.
    \end{cases}
\end{align*}
\end{theorem}

\begin{proof}

We first prove the interior equation. By Corollary \ref{Cor:Green:l2diff}, there exists
\(w'\in L^2(D_2)\) such that
$$
\frac{\tilde w^h-\tilde w^0}{h}\to w'
\quad\text{in }L^2(D_2).
$$
Let $\psi\in C_c^\infty(\Omega^0)$ and set
$K=\operatorname{supp}\psi$. Since $K\Subset\Omega^0$, for all
sufficiently small $|h|$ we have $K\subset\Omega^h$. Therefore,
$w^h$ and $w^0$ are harmonic in a neighborhood of $K$, and hence
$$
q_h:=\frac{w^h-w^0}{h}
$$
satisfies $\Delta q_h=0$ in $K$. Thus
$$
\int_K q_h\Delta\psi\,dx=0.
$$
Passing to the limit using $q_h\to w'$ in $L^2(K)$, we obtain
$$
\int_{\Omega^0}w'\Delta\psi\,dx=0.
$$
Since \(\psi\in C_c^\infty(\Omega^0)\) was arbitrary, it follows that
$$
\Delta w'=0\quad\text{in }\Omega^0,
$$
by Weyl's lemma. Since
$G^t_\lambda=w^t+\Phi$ and $\Phi$ is independent of $t$, we have
$G_\lambda'=w'$. Hence,
$$
\Delta G_\lambda'=0\quad\text{in }\Omega^0.
$$

We next prove the boundary condition on $\Gamma_S$. Since
$\operatorname{supp}V\Subset D_2$, there exists a neighborhood $U$ of
$\Gamma_S$ such that $V=0$ in $U$. Hence $T_t=\operatorname{Id}$
in $U$, and $U\cap\Omega^t=U\cap\Omega^0$ for all sufficiently small
$|t|$.

Let $\psi\in H^1(\Omega^0)$ with $\operatorname{supp}\psi\subset U$.
Then $\psi$ is an admissible test function for the weak formulation of
$w^t$ for all small $|t|$. Therefore
$$
\int_{\Omega^0}\nabla w^t\cdot\nabla\psi\,dx
=
\int_{\Gamma_S}
(\lambda w^t+\lambda\Phi-\partial_\nu\Phi)\psi\,dS.
$$
Subtracting the equation at $t=0$ and dividing by $t$, we get
$$
\int_{\Omega^0}
\nabla\left(\frac{w^t-w^0}{t}\right)\cdot\nabla\psi\,dx
=
\lambda\int_{\Gamma_S}
\frac{w^t-w^0}{t}\psi\,dS.
$$
Since $T_t=\operatorname{Id}$ in $U$, the $H^1$-differentiability
of $w^t\circ T_t$ implies
$$
\frac{w^t-w^0}{t}\to w'
\quad\text{in }H^1(U\cap\Omega^0).
$$
Passing to the limit gives
$$
\int_{\Omega^0}\nabla w'\cdot\nabla\psi\,dx
=
\lambda\int_{\Gamma_S}w'\psi\,dS.
$$
Since \(G_\lambda'=w'\), we have
$$
\int_{\Omega^0}\nabla G_\lambda'\cdot\nabla\psi\,dx
=
\lambda\int_{\Gamma_S}G_\lambda'\psi\,dS.
$$
Together with $\Delta G_\lambda'=0$ in $\Omega^0$, Green's identity yields
$$
\int_{\Gamma_S}(\partial_\nu G_\lambda'-\lambda G_\lambda')\psi\,dS=0.
$$
Since the trace of $\psi$ on $\Gamma_S$ is arbitrary, it follows that
$$
\partial_\nu G_\lambda'=\lambda G_\lambda'
\quad\text{on }\Gamma_S.
$$

We then prove the the Dirichlet boundary condition on $\Gamma_D$. From now on, we assume $x$ is in the small neighborhood of $\Gamma_D$. By the chain rule, we have
\begin{align*}
     (\bar{w}^0)'=\frac{d}{dt}\tilde{w}^t(T(t,x))\Bigl |_{t=0}&=(\tilde{w}^{0})'(x) + (\nabla \tilde{w}^{0}(x) \cdot d)\\
    &=(\tilde{w}^{0})'(x) + (\nabla G^{0}_\lambda(x) \cdot d) -\nabla\Phi(x)\cdot d.
\end{align*}

Note that $\bar{w}^t$ is differentiable in $H^1(D_2)$ (see Corollary~\ref{Cor:green:H1diff}). So we have
$$(\bar{w}^0)'=(\tilde{w}^{0})'(x) + (\nabla G^{0}_\lambda(x) \cdot d)-\nabla\Phi(x)\cdot d\in H^1(D_2).$$
From proof of Corollary \ref{Cor:green:H1diff}, trace of $(\bar{w}^0)'$ on $\Gamma_D^0$ is $-\nabla\Phi(x)\cdot d$. Therefore, we have
$$G_\lambda'(x)=(\tilde{w}^{0})'(x) =- (\nabla G^{0}_\lambda(x) \cdot d) \t{ for } x\in \Gamma_D^0$$
This shows $G_\lambda'$ satisfies the desired Dirichlet condition on $\Gamma_D$.
\end{proof}


\begin{theorem}\label{Thm:Green:Jdiff}
The functional $t \mapsto \int_{\Gamma_S}\left(G_\lambda^t(x)\right)^2\,dS$ is differentiable at $t=0$ and the derivative is $2\int_{\Gamma_S} G_\lambda^0 G_\lambda'\,dS$.
\end{theorem}
\begin{proof}
Let $$A: H^1(D_2)\to \R, \ A(f)= \int_{\Gamma_S} (f)^2\,dS \ \t{ and } \ B: H^1(D_2)\to \R, \ B(f)= \int_{\Gamma_S} f\Phi\, dS.$$ Then $A$ and $B$ are differentiable, so $(A+2B)(\bar{w}^t)+\int_{\Gamma_S} \Phi^2\,dS: \R \to \R$ is also differentiable. Using the chain rule, we have $(\int_{\Gamma_S} (G_\lambda^t)^2\,dS)'\big |_{t=0}=2\int_{\Gamma_S} G_\lambda^0G_\lambda'\, dS$.
\end{proof}

\section{Shape optimization for eigenvalues} \label{sec:optm}

We now develop a theoretical framework for shape optimization of the domain configuration for maximizing a simple eigenvalue.
For the detailed description, we write $D_1^q$ for the translate of $D_1$ by a vector $q \in \mathbb{R}^n$, and set
$$
    \Omega^q := D_2 \setminus \overline{D_1^q}.
$$
For each admissible translation \(q\), we define
\[
    J_\lambda(q) := \int_{\Gamma_S} \left(G_\lambda^q(x_s,x)\right)^2 \, dS.
\]
Theorem~\ref{thm:Numerical:J_expression} implies that, if $\mu$ is an eigenvalue of multiplicity $m$, then
$$
\lim_{\lambda\to\mu}
(\lambda-\mu)^2 J_\lambda
=
\sum_{\sigma_j=\mu}|u_j(x_s)|^2.
$$

Applying Theorem~\ref{Thm:Green:Jdiff} to translations of the inner boundary in the $q_i$-direction, we obtain
$$
\frac{\partial J_\lambda}{\partial q_i}(q)
=
2\int_{\Gamma_S}G_\lambda^q(x)
\frac{\partial G_\lambda^q}{\partial q_i}(x),dS,
\qquad i=1,\ldots,n.
$$
The following lemma establishes differentiability of an eigenfunction evaluated at a fixed interior point. From this lemma, we have remainder term
\begin{align} \label{def:optm:R_lambda}
    R_\lambda(q)
    :=
    J_\lambda(q)
    -
    \frac{|u_k^q(x_s)|^2}{(\sigma_k^q-\lambda)^2}
\end{align}
is differentiable.

\begin{lemma}
\label{lem:interior-point-differentiability}
Suppose that there exists $r>0$, $\varepsilon>0$ and $x_s$
such that
\[
    \overline{B_r(x_s)}\subset\Omega^t,\quad  x_s \in \Omega^t \setminus \t{supp} (V)
    \qquad\text{for every }t\in (-\varepsilon,\varepsilon).
\]
Then for simple eigenfunction $u_k$, the function $a(t):=u_k^t(x_s)$
is a $C^1$ function.
\end{lemma}

\begin{proof}
Write $B:=B_r(x_s)$ and define the averaging functional
\[
    \mathcal M:L^2(B)\longrightarrow\mathbb R,
    \qquad
    \mathcal M(v):=\frac{1}{|B|}\int_B v(y)\,dy.
\]
By the Cauchy--Schwarz inequality,
\[
    |\mathcal M(v)|
    \le |B|^{-1/2}\|v\|_{L^2(B)}.
\]
Thus, $\mathcal M$ is bounded and linear.

Let $\widehat{u}_k^t:=u_k^t\circ T_t$. Then for sufficiently small $r>0$, we have $u^t_k|_B=\widehat{u}_k^t|_B$.
Since $u^t_k$ is harmonic in a neighborhood of $\overline B$,
the mean value property gives
\[
    a(t)=u_k^t(x_s)
    =\frac{1}{|B|}\int_B \widehat{u}_k^t(y)\,dy
    =\mathcal M( u_k^t).
\]
Consequently, $a(t)=\mathcal M\circ \widehat{u}_k^t$ is a $C^1$ function, since $\widehat{u}_k^t$ is a $C^1$ function by Theorem \ref{thm:prop:diff}.
\end{proof}

\begin{theorem}\label{thm:gradient-alignment}
For a simple eigen pair $\sigma_k^{q_0}$ and $u_k^{q_0}$, assume
\[
    u_k^{q_0}(x_s)\neq 0,
    \qquad
    \nabla_q\sigma_k^{q_0}\neq 0.
\]
Set $\mu_0:=\sigma_k^{q_0}$ and
assume that there exists $\varepsilon>0$ such that
$(\mu_0-\varepsilon,\mu_0+\varepsilon)$ contains no other Steklov--Dirichlet
eigenvalue of $\Omega^{q_0}$ and satisfies
$$\sup_{0<|\lambda-\mu_0|<\varepsilon}\lVert \nabla_q R_\lambda(q_0)\rVert < \infty.$$
Then, as $\lambda\to \sigma_k^{q_0}$,
\[
    \nabla_qJ_\lambda(q_0)
    =
    \frac{2|u_k^{q_0}(x_s)|^2}
    {(\lambda-\sigma_k^{q_0})^3}
    \nabla_q\sigma_k^{q_0}
    +
    O\left(
    (\lambda-\sigma_k^{q_0})^{-2}
    \right).
\]
Consequently,
\begin{align*}
    \frac{\nabla_qJ_\lambda(q_0)}
    {|\nabla_qJ_\lambda(q_0)|}
    \longrightarrow
    \pm\frac{\nabla_q\sigma_k^{q_0}}
    {|\nabla_q\sigma_k^{q_0}|}
    \qquad
    \text{as }
    \lambda\to (\sigma_k^{q_0})^\pm.
\end{align*}

\end{theorem}

\begin{proof}
Set
\[
    a_k(q):=|u_k^q(x_s)|^2
\]
and
\[
    \delta:=\lambda-\sigma_k^{q_0}.
\]
By the spectral representation (see also \eqref{def:optm:R_lambda}), 
\[
    J_\lambda(q)
    =
    \frac{a_k(q)}
    {(\sigma_k^q-\lambda)^2}
    +
    R_\lambda(q).
\]
We obtain
\begin{align*}
    \nabla_qJ_\lambda(q)
    &=
    \frac{\nabla_q a_k(q)}
    {(\sigma_k^q-\lambda)^2}
    -
    \frac{2a_k(q)}
    {(\sigma_k^q-\lambda)^3}
    \nabla_q\sigma_k^q
    +
    \nabla_qR_\lambda(q).
\end{align*}
Evaluating at $q=q_0$ and using
\[
    \sigma_k^{q_0}-\lambda=-\delta,
\]
we obtain,
\begin{align*}
    \nabla_qJ_\lambda(q_0)
    &=
    \frac{\nabla_q a_k(q_0)}{\delta^2}
    +
    \frac{2a_k(q_0)}{\delta^3}
    \nabla_q\sigma_k^{q_0}
    +
    \nabla_qR_\lambda(q_0),
\end{align*}
where differentiability of $a_k$ is guaranteed by Lemma \ref{lem:interior-point-differentiability}.
Since $a_k$ is $C^1$ and
$\nabla_qR_\lambda(q_0)$ remains bounded as
$\lambda\to \sigma_k^{q_0}$, it follows that
\[
    \nabla_qJ_\lambda(q_0)
    =
    \frac{2a_k(q_0)}{\delta^3}
    \nabla_q\sigma_k^{q_0}
    +
    O(\delta^{-2}).
\]

We may factor the leading term as
\[
    \nabla_qJ_\lambda(q_0)
    =
    \frac{2a_k(q_0)}{\delta^3}
    \left(
        \nabla_q\sigma_k^{q_0}
        +O(\delta)
    \right).
\]
Because $\nabla_q\sigma_k^{q_0}\neq0$, normalization yields
\[
    \frac{\nabla_qJ_\lambda(q_0)}
    {|\nabla_qJ_\lambda(q_0)|}
    \longrightarrow
    \pm\frac{\nabla_q\sigma_k^{q_0}}
    {|\nabla_q\sigma_k^{q_0}|}
\]
as $\delta\to 0^\pm$, which proves the result.
\end{proof}

We now provide an error estimate for the eigenvalue estimator based on finite differences of $H(\lambda):= \sqrt{J_\lambda}$.

\begin{theorem}
\label{thm:pole-estimator}
Fix an admissible translation $q$, and let $\mu=\sigma_k^q$
be a simple Steklov--Dirichlet eigenvalue with
$u_k^q(x_s)\neq0$.
Then there exist constants $C>0$ and $\varepsilon>0$ such
that, for every real $\lambda$ and $h>0$ satisfying
\[
    \delta:=|\lambda-\mu|>0,
    \qquad
    \delta+h<\varepsilon,
    \qquad
    (\lambda-\mu)(\lambda+h-\mu)>0,
\]
the estimator
\[
    \widehat{\mu}_h
    :=
    \lambda-
    \frac{h}{H(\lambda)/H(\lambda+h)-1}
\]
is well defined and satisfies
\[
    |\widehat{\mu}_h-\mu|
    \le C\,\delta(\delta+h)^2.
\]
\end{theorem}

\begin{proof}
Throughout the proof, $q$ is fixed.
Let $a:=|u_k^q(x_s)|^2>0$.
By Theorem~\ref{thm:Numerical:J_expression},
\[
    J_\lambda
    =
    \frac{a}{(\lambda-\mu)^2}+R(\lambda),
    \qquad
    R(\lambda)
    :=
    \sum_{j\ne k}
    \frac{|u_j(x_s)|^2}{(\sigma_j-\lambda)^2}.
\]
Since $\mu$ is isolated, the series defining $R$ and
its first $\lambda$-derivative converge uniformly on a
sufficiently small closed interval about $\mu$.
Thus, $R$ is a $C^1$ function near $\mu$.

For $s$ near zero, define
\[
    \mathcal F(s)
    :=
    \frac{s}{\sqrt a}
    \left(
        1+\frac{s^2}{a}R(\mu+s)
    \right)^{-1/2}.
\]
The expression inside the parentheses is positive for
sufficiently small $|s|$.
Moreover, $\mathcal F$ is $C^1$ near zero and satisfies
\[
    \mathcal F(s)
    =
    \frac{s}{\sqrt a}\bigl(1+O(s^2)\bigr),
    \qquad
    \mathcal F'(s)
    =
    \frac{1}{\sqrt a}\bigl(1+O(s^2)\bigr).
\]
In particular, there exists $\varepsilon$ such that
\[
    \mathcal F'(s)\ge\frac{1}{2\sqrt a}
    \qquad\text{for }|s|<\varepsilon.
\]
For $s\ne0$, we also have
\[
    |\mathcal F(s)|=\frac{1}{H(\mu+s)}.
\]

Set $s_0:=\lambda-\mu$, so that $\delta=|s_0|$.
Since $s_0$ and $s_0+h$ have the same sign,
\[
    \frac{\mathcal F(s_0+h)}{\mathcal F(s_0)}
    =
    \frac{H(\lambda)}{H(\lambda+h)}.
\]
By the mean value theorem, there exists
$\xi\in(s_0,s_0+h)$ such that
\[
    \mathcal F(s_0+h)-\mathcal F(s_0)
    =
    h\mathcal F'(\xi).
\]
The lower bound on $\mathcal F'$ shows that this difference
is nonzero. Hence the estimator is well defined.

Writing $\widehat s_h:=\lambda-\widehat\mu_h$, we obtain
\[
\begin{aligned}
    \widehat s_h
    &=
    h\frac{\mathcal F(s_0)}
    {\mathcal F(s_0+h)-\mathcal F(s_0)}\\
    &=
    \frac{\mathcal F(s_0)}{\mathcal F'(\xi)}.
\end{aligned}
\]
Since
\[
    |\xi|
    \le \max\{|s_0|,|s_0+h|\}
    \le \delta+h,
\]
the preceding estimates yield
\[
    \widehat s_h
    =
    s_0\,
    \frac{1+O(\delta^2)}
         {1+O((\delta+h)^2)}
    =
    s_0\bigl(1+O((\delta+h)^2)\bigr).
\]
Consequently,
\[
    |\widehat\mu_h-\mu|
    =
    |s_0-\widehat s_h|
    \le C\,|s_0|(\delta+h)^2
    =
    C\,\delta(\delta+h)^2,
\]
which proves the result.
\end{proof}

\FloatBarrier
\section{Numerical computation of eigenvalues} \label{sec:num}
In this section, we present a numerical algorithm to determine the location of the inner boundary of a shell-type domain that maximizes the Steklov–Dirichlet eigenvalue. In Section \ref{subsec:num:framework}, we describe the computational framework for evaluating the Steklov–Dirichlet eigenvalue and the Green function introduced in Section \ref{sec:spectral}, using the layer potential approach. In Section \ref{subsec:num:algorithm}, we present the algorithm, based on the results of Sections \ref{sec:derv} and \ref{sec:optm}. Finally, we provide numerical examples in Section \ref{subsec:num:ex}.

\subsection{Boundary Integral Formulation} \label{subsec:num:framework}
We present the layer potential approach to compute the Steklov--Dirichlet eigenvalues and Green functions. See \cite{Colton:1998:IAE} for a general introduction to the layer potential technique, and \cite[Section~6]{Ammari:2020:OSE} for its application to the case of the mixed Steklov--Neumann problem.

The single-layer potential is defined by
$$\mathcal{S}[\varphi](x):=\int_{\partial\Omega} \Phi(x,y)\varphi(y) \, d\sigma_y, \quad  x\in\Omega,$$
where $\varphi\in L^2(\partial \Omega)$ is a density function. The normal derivative of the single-layer potential satisfies the jump relation, given by
$$\partial_\nu \mathcal{S}[\varphi](x) = (K^*-\frac12 I)[\varphi](x)=\int_{\partial\Omega}\partial_{\nu_x}\Phi(x,y)\varphi(y)d\sigma_y -\frac{1}{2}\varphi(x), \quad  x\in\p \Omega.$$
A harmonic function in $\Omega$ can be represented as a single layer potential with an appropriate density function.

In what follows, we decompose the density function into two components associated with two boundary components of $\p \Omega$.
Let \(\varphi_1\) and \(\varphi_2\) denote the density components for $\Gamma_D$ and $\Gamma_S$, respectively. Then we obtain the boundary integral equation
\[
\begin{bmatrix}
 P_{1}\mathcal{S} & P_{1}\mathcal{S}\\
P_{2}K^* & P_{2}(K^*-\frac12 I)
\end{bmatrix}
\begin{bmatrix}
\varphi_1\\ \varphi_2
\end{bmatrix}
=
\sigma
\begin{bmatrix}
 0 & 0\\
P_{2}\mathcal{S} & P_{2}\mathcal{S}
\end{bmatrix}
\begin{bmatrix}
\varphi_1\\ \varphi_2
\end{bmatrix}
\]
where $P_1,P_2$ are projections onto $\Gamma_D, \Gamma_S$, respectively. 
After discretization, this equation yields a generalized eigenvalue problem, which can be solved using MATLAB's
\texttt{eig} function.

We also express the Green functions using single layer potentials. The key idea is to separate the singularity. Let $\psi_1$ and $\psi_2$ denote the density components of $G^t_{\lambda}-\Phi$ supported on the respective boundary. Then we have
\begin{equation*}
    \begin{bmatrix} 
   P_{1}\mathcal{S} & P_{1}\mathcal{S} \\
   P_{2}(K^* -\lambda\mathcal{S}) & P_{2}(K^*-\frac12 I -\lambda\mathcal{S})
   \end{bmatrix}
   \begin{bmatrix} 
   \psi_1 \\
   \psi_2
   \end{bmatrix}
   =
   \begin{bmatrix} 
   -P_{1}\Phi \\
   P_{2}(\lambda \Phi-\frac{\partial \Phi}{\partial \nu}).
   \end{bmatrix}
\end{equation*}  
Since $\lambda$ is fixed, computation of the density reduces to solve the boundary integral equation. 

Additionally, recall that the shape derivative of the Green function satisfies the following equation (see Theorem \ref{thm:derv:G'}),
\begin{align*}
    \begin{cases}
        \Delta G_\lambda'(x) = 0 & \text{in } \Omega^{0},\\
  G_\lambda'(x) = - \frac{\partial G_\lambda^0}{\partial \nu} (\nu \cdot d)       & \text{on } \Gamma_D^{0},\\[2mm]
  \dfrac{\partial G_\lambda'}{\partial \nu}(x) = \lambda G_\lambda'
                 & \text{on } \Gamma_S.
    \end{cases}
\end{align*}
The shape derivative of the Green function is computed by the equation
\begin{equation*}
    \begin{bmatrix}
        P_{1}\mathcal{S} & P_{1}\mathcal{S}\\
        P_{2}(K^* -\lambda\mathcal{S})
        &
        P_{2}(K^*-\frac12 I -\lambda\mathcal{S})
    \end{bmatrix}
    \begin{bmatrix}
        \zeta_1\\
        \zeta_2
    \end{bmatrix}
    =
    \begin{bmatrix}
        -\partial_\nu G_\lambda^0 \,(d\cdot\nu)|_{\Gamma_D^0}\\
        0
    \end{bmatrix},
\end{equation*}
for a translation of inner boundary in the direction $d$.

\subsection{Algorithm} \label{subsec:num:algorithm}

The algorithm moves the inner boundary along the steepest descent direction of the functional $J_\lambda$, as summarized in Algorithm \ref{algm}. 
Since Theorem \ref{thm:gradient-alignment} states that $$\frac{\nabla_q J_\lambda(q)}{\norm{\nabla_q J_\lambda(q)}} \rightarrow -\frac{\nabla_q \sigma_1(q)}{\norm{\nabla_q \sigma_1(q)}} \quad 
\text{ as } \lambda \to (\sigma_1)^-,$$ 
decreasing $J_\lambda$ has the effect of driving $\sigma_1^{q}$ to the gradient direction if $\lambda < \sigma_1$.
For the first iteration, we compute $\sigma_1^{q_0}$ and set $\lambda^{q_0}=0.99 \, \sigma_1^{q_0}$. 
The algorithm then numerically yields an increase in the first eigenvalue.
As demonstrated by the numerical simulations in Section \ref{subsec:num:ex}, it drives the first eigenvalue toward its approximated maximum. 

Specifically, in the Algorithm \ref{algm}, we compute gradient descent direction $d_\ell$ using Theorem \ref{Thm:Green:Jdiff} and estimated the first eigenvalue using Theorem \ref{thm:pole-estimator}.
In our implementation, the step size $\eta>0$ is chosen sufficiently small relative to the size of $D_2$, the region enclosed by the outer boundary $\Gamma_S$. For instance, one may set $\eta$ to one tenth of the radius of the largest inscribed circle of $D_2$.



\begin{algorithm}
\caption{Translation search for maximizing the first
Steklov--Dirichlet eigenvalue}
\label{algm}
\begin{algorithmic}[1]

\Require Initial translation $q_0\in\mathcal A$, finite-difference
parameter $h>0$, and step size $\eta>0$.

\State Compute $\sigma_1^{q_0}$ and set
\[
    \lambda_0=0.99\,\sigma_1^{q_0},
    \qquad d_{-1}=0.
\]

\For{$\ell\geq 0$}

    \State Compute the translation gradient $g_\ell=\nabla_q J_{\lambda_\ell}(q_\ell)$ (see Theorem  \ref{Thm:Green:Jdiff}).

    \State Set $d_\ell=\frac{-g_\ell}{\|g_\ell\|}.$

    \If{$d_{\ell-1}\cdot d_\ell<0$}
        \State Set step size $\eta = \eta/2$. 
    \EndIf

    \State Update $q_{\ell+1}=q_\ell+\eta d_\ell$.
    \If{$\eta < 0.05$}
        \State \Return $q_m=q_{\ell+1}$.
    \EndIf
    \State Compute estimated eigenvalue $\widehat{\mu}_h^{q_{\ell+1}}(\lambda_\ell)$ (see Theorem  \ref{thm:pole-estimator}).
    \State Update $\lambda_{\ell+1}=0.99\widehat{\mu}_h^{q_{\ell+1}}(\lambda_\ell)$.
    \State Update $\ell=\ell+1$.
\EndFor

\end{algorithmic}
\end{algorithm}

\FloatBarrier
\subsection{Examples}  \label{subsec:num:ex}
We now perform numerical simulations to find the location of the inner boundary $\Gamma_D$ that maximizes the first eigenvalue on $\Omega$.
We consider two-dimensional domains
$\Omega^q=D_2\setminus\overline{D_1^q}$ for a vector $q \in \R^2$ (see Section \ref{sec:optm}.)
We employ the layer-potential formulation described in Section \ref{subsec:num:framework}.

The finite-difference parameter $h$ in Algorithm~\ref{algm} is fixed by
$h=10^{-4}$ in the simulation. The inner and outer boundaries are discretized using $N_D=300$ and $N_S=1500$ nodes, respectively, with uniform sampling. Adjacent nodes are connected by straight segments to form closed polygonal approximations of the boundaries. The initial translation step size is $\eta=0.2$. The source point $x_s$ is chosen to remain in $\Omega^q$ throughout the iteration. The shape and orientation of the inner boundary are fixed throughout the translation search. In the figures, we represent the trajectory along Algorithm~\ref{algm}. The gray inner boundary indicates the initial configuration, while the black inner boundary represents the final configuration.

\subsubsection{Elliptical outer boundary}

\begin{figure}[h]
  \centering

  \subfloat[Circular inner boundary]{%
    \begin{minipage}{0.4\textwidth}
      \centering
      \includegraphics[width=\linewidth]{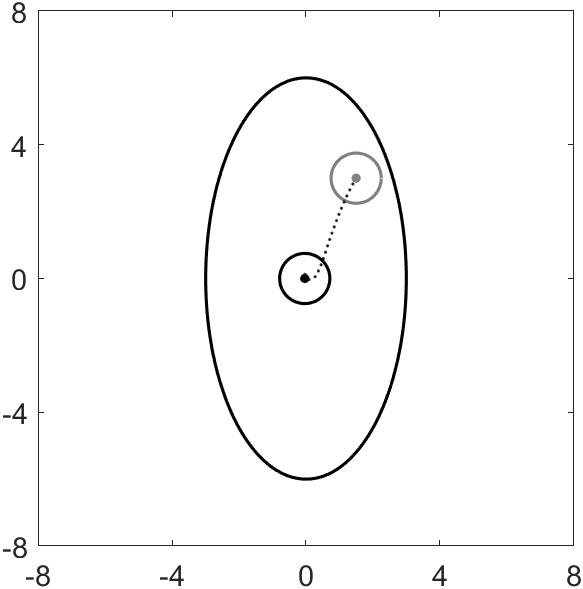}
    \end{minipage}
  }
  \hfill
  \subfloat[Elliptical inner boundary]{%
    \begin{minipage}{0.4\textwidth}
      \centering
      \includegraphics[width=\linewidth]{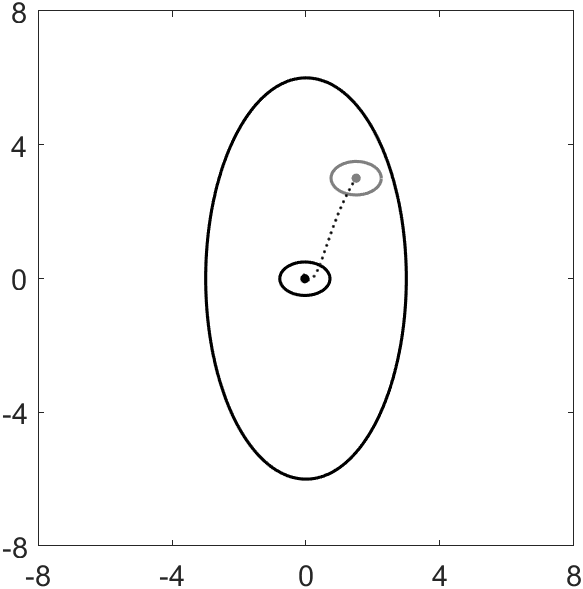}
    \end{minipage}
  }

  \caption{Translation search with an elliptical outer boundary.}
  \label{fig:ex:ellipse}
\end{figure}

We first consider an elliptical outer boundary and compare two fixed
inner shapes: a circle and an ellipse.
Figure~\ref{fig:ex:ellipse} shows the trajectory of the inner boundary.
Table~\ref{tab:ellipse} report
selected iterations for the two configurations. We fix $x_s=(-2.5,1)$ for each case. 
The outer boundary is an ellipse with semiaxes $3$ and $6$, while inner boundaries are a circle with radius $0.75$ and an ellipse with semiaxes $0.5$ and $0.75$, respectively.

In both cases, Algorithm~\ref{algm} moves the inner boundary toward the center of the outer ellipse, following similar trajectories. The estimated first eigenvalues for the final configurations are $0.1099$ and $0.1008$ for the two cases, respectively, both of which agree with the reference values to four decimal places. 

\FloatBarrier
\begin{table}[h]
  \centering

  \subfloat[Circular inner boundary]{%
    \begin{minipage}{0.47\textwidth}
      \centering
      \begin{tabular}{cccc}
        \hline
        $\ell$ & $q_x$ & $q_y$ & $\sigma_1^{q_\ell}$ \\
        \hline
        0  & 1.50 & 3.00 & 0.0592 \\
        1  & 1.40 & 2.83 & 0.0618 \\
        5  & 1.06 & 2.10 & 0.0731 \\
        10 & 0.71 & 1.17 & 0.0898 \\
        15 & 0.37 & 0.23 & 0.1075 \\
        20 & -0.04 & -0.08 & 0.1097 \\
        23 & \bf{-0.04} & \bf{0.00} & $\mathbf{0.1099}$ \\
        \hline
      \end{tabular}
    \end{minipage}
  }
  \hfill
  \subfloat[Elliptical inner boundary]{%
    \begin{minipage}{0.47\textwidth}
      \centering
      \begin{tabular}{cccc}
        \hline
        $\ell$ & $q_x$ & $q_y$ & $\sigma_1^{q_\ell}$ \\
        \hline
        0  & 1.50  & 3.00 & 0.0562 \\
        1  & 1.40  & 2.82 & 0.0587 \\
        5  & 1.04  & 2.11 & 0.0692 \\
        10 & 0.69  & 1.18 & 0.0843 \\
        15 & 0.34  & 0.24 & 0.0990 \\
        20 & 0.01  & -0.08 & 0.1007 \\
        23 & \bf{-0.03} & \bf{-0.01} & $\mathbf{0.1008}$ \\
        \hline
      \end{tabular}
    \end{minipage}
  }

  \caption{Selected iterates of the translation search.}

  \label{tab:ellipse}
\end{table}


\subsubsection{Kite-shaped outer boundary} 

\begin{figure}[h]
  \centering
  \subfloat[Symmetric kite shape]{%
    \includegraphics[width=0.4\textwidth]{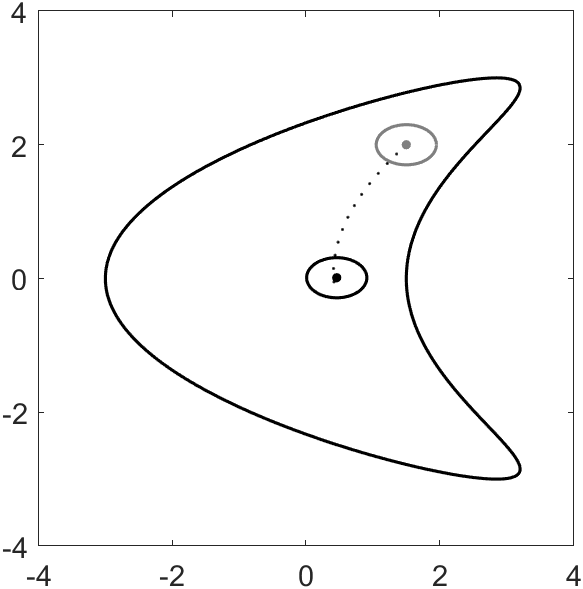}
  }
  \hfill
  \subfloat[Asymmetric kite shape]{%
    \includegraphics[width=0.4\textwidth]{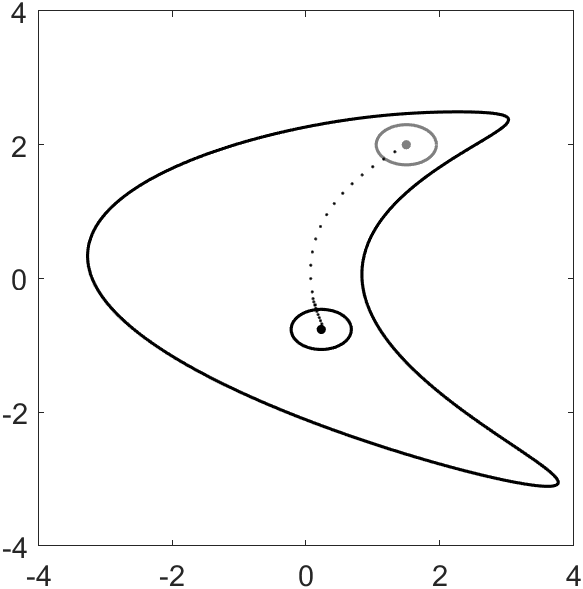}
  }
  \caption{Translation paths for a kite-shaped outer boundary.}
  \label{fig:kite-paths}
\end{figure}

We next consider kite-shaped outer boundaries with elliptical inner boundary.
Figure~\ref{fig:kite-paths} presents the trajectory along Algorithm~\ref{algm}. The left panel shows a kite-shaped boundary symmetric about the \(x\)-axis, whereas the right panel shows an asymmetric boundary. 
We use
$$r\left[ 0.35 + 0.75\cos\theta - 0.6\cos(2\theta),\sin\theta \right], \quad 0\le\theta<2\pi,$$
and
$$r\left[ 0.35 + 0.75\cos\theta - 0.7\cos(2\theta) + 0.24\sin(3\theta),0.1+\sin\theta+0.21\cos(2\theta) \right], \quad 0\le\theta<2\pi,$$
respectively, with $r=2.8$. For both shapes, the inner boundary is an ellipse with semiaxes $0.45$ and $0.3$. 
We choose the fixed source point $x_s=(-2,0)$. 

In both cases, Algorithm~\ref{algm} places the inner boundary in the right-hand part of the kite-shaped outer domain. In the asymmetric case, the final position lies below that obtained in the symmetric case. For both kite-shaped outer boundaries, the candidate configurations for the maximum are shifted to the right. The final vertical displacement is lower in the asymmetric case than in the symmetric case. 

Notably, for the asymmetric kite shape, we obtain $\lambda_{15}\approx 0.0759$. This value exceeds the first eigenvalue $\sigma_1^{q_{15}}=0.0756$, violating the intended condition that the spectral parameter remain below the first eigenvalue. At the next iteration, the computed search direction deviated from the intended ascent direction for the first eigenvalue. This temporary deviation occurred only once in this experiment, and the expected search direction was recovered in subsequent iterations.
This observation illustrates a limitation of the spectral parameter update: multiplying the estimated eigenvalue by $0.99$ may not ensure that the updated parameter remains below the true eigenvalue. If the estimation error exceeds the margin provided by this factor, the parameter may cross the first eigenvalue, changing the sign of the leading term in the gradient alignment formula.

\begin{table}[h]
  \centering
  \subfloat[Symmetric kite shape]{%
    \begin{minipage}{0.47\textwidth}
      \centering
      \begin{tabular}{cccc}
        \hline
        $\ell$ & $q_x$ & $q_y$ & $\sigma_1^{q_\ell}$ \\
        \hline
        0  & 1.50 & 2.00 & 0.0546 \\
        1  & 1.36 & 1.86 & 0.0592 \\
        5  & 0.83 & 1.26 & 0.0802 \\
        10 & 0.44 & 0.35 & 0.1097 \\
        15 & \bf{0.46} & \bf{0.01} & $\mathbf{0.1162}$ \\
        \hline
      \end{tabular}
    \end{minipage}
  }
  \hfill
  \subfloat[asymmetric kite shape]{%
    \begin{minipage}{0.47\textwidth}
      \centering
      \begin{tabular}{cccc}
        \hline
        $\ell$ & $q_x$ & $q_y$ & $\sigma_1^{q_\ell}$ \\
        \hline
        0  & 1.50  & 2.00 & 0.0333 \\
        1  & 1.33  & 1.89 & 0.0359 \\
        5  & 0.69  & 1.55 & 0.0466 \\
        10 & 0.15  & 0.59 & 0.0597 \\
        15 & 0.14  & -0.39 & 0.0756 \\
        20 & 0.16  & -0.49 & 0.0776 \\
        25 & 0.25  & -0.73 & 0.0821 \\
        27 & \bf{0.23} & \bf{-0.76} & $\mathbf{0.0822}$ \\
        \hline
      \end{tabular}
    \end{minipage}
  }

  \caption{Selected iterates for the translation searches.}
  \label{tab:kite-results}
\end{table}

\subsubsection{Other boundaries} 

\begin{figure}[h]
  \centering

  \subfloat[Butterfly shape]{%
    \includegraphics[width=0.4\textwidth]{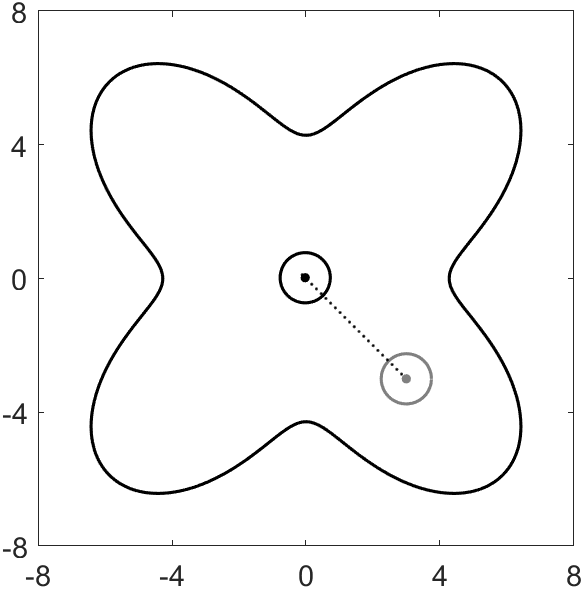}
  }
  \hfill
  \subfloat[Dumbbell shape]{%
    \includegraphics[width=0.4\textwidth]{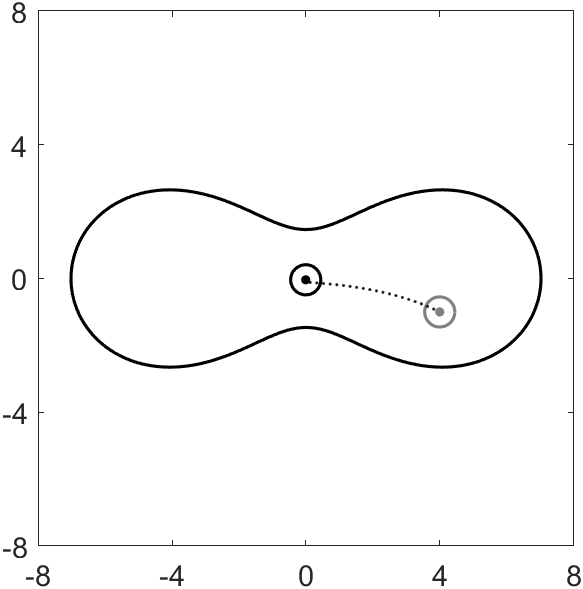}
  }

  \caption{Translation paths for various outer boundaries.}
  \label{fig:other-paths}
\end{figure}

Lastly, we consider butterfly-shaped and dumbbell-shaped outer boundaries.
Figure~\ref{fig:other-paths} presents the corresponding translation paths, and Table~\ref{tab:other-results} lists selected iteration histories. The radii of the inner boundaries are $0.5$ and $0.3$. We fix $x_s=(-4,4)$ for the butterfly case and $x_s=(-6,0)$ for the dumbbell case. 

For both the butterfly-shaped and dumbbell-shaped outer boundaries, the numerical results indicate that the first eigenvalue is maximized when the inner boundary is located near the center. This observation may be unexpected, as one might anticipate the maximum to occur when the inner boundary is placed within one of the lobes. Thus, despite the non-convexity of the outer boundaries and their distinct geometric features, the maximizing configuration remains near the center in both cases.

\begin{table}[h]
  \centering

  \subfloat[Butterfly shape]{%
    \begin{minipage}{0.47\textwidth}
      \centering
      \begin{tabular}{cccc}
        \hline
        $\ell$ & $q_x$ & $q_y$ & $\sigma_1^{q_\ell}$ \\
        \hline
        0  & 3.00 & -3.00 & 0.0308 \\
        1  & 2.86 & -2.86 & 0.0320 \\
        5  & 2.29 & -2.29 & 0.0371 \\
        10 & 1.59 & -1.59 & 0.0439 \\
        15 & 0.88 & -0.88 & 0.0499 \\
        20 & 0.17 & -0.17 & 0.0534 \\
        25 & -0.01 & 0.01 & 0.0536 \\
        27 & \bf{-0.02} & \bf{0.02} & $\mathbf{0.0536}$ \\
        \hline
      \end{tabular}
    \end{minipage}
  }
  \hfill
  \subfloat[Dumbbell shape]{%
    \begin{minipage}{0.47\textwidth}
      \centering
      \begin{tabular}{cccc}
        \hline
        $\ell$ & $q_x$ & $q_y$ & $\sigma_1^{q_\ell}$ \\
        \hline
        0  & 4.00  & -1.00 & 0.0279 \\
        1  & 3.82  & -0.91 & 0.0286 \\
        5  & 3.07  & -0.63 & 0.0320 \\
        10 & 2.11  & -0.37 & 0.0379 \\
        15 & 1.12  & -0.21 & 0.0473 \\
        20 & 0.12  & -0.12 & 0.0612\\
        24 & \bf{-0.01} & \bf{-0.04} & $\mathbf{0.0620}$ \\
        \hline
      \end{tabular}
    \end{minipage}
  }

  \caption{Selected iterates for the translation searches.}
  \label{tab:other-results}
\end{table}


\section{Conclusion} \label{sec:concl}
In this paper, we developed a Green function-based approach to the spectral analysis and optimization of Steklov–Dirichlet eigenvalues on smooth shell-type domains. We established a spectral representation of the associated Green function and derived its shape derivative under translations of the inner boundary. Based on these results, we proposed an algorithm to determine the position of the inner boundary that maximizes the first eigenvalue without repeatedly solving the eigenvalue problem. In particular, we derived an eigenvalue estimator with a local error bound to approximate the first eigenvalue during the optimization process. As a future direction, we will provide a convergence analysis of the proposed algorithm. 


\bibliographystyle{amsplain} 
\bibliography{2026_SDEigenvalue_genShape}

@article {Agranovich:2006:MPS,
    AUTHOR = {Agranovich, Mikhail Semenovich},
     TITLE = {On a mixed {P}oincar\'{e}-{S}teklov type spectral problem in a
              {L}ipschitz domain},
   JOURNAL = {Russ. J. Math. Phys.},
  FJOURNAL = {Russian Journal of Mathematical PhysiDcs},
    VOLUME = {13},
      YEAR = {2006},

     PAGES = {239--244},
      ISSN = {1061-9208},
   MRCLASS = {35J25 (35P20 47F05)},
    DOI = {10.1134/S1061920806030010},
}

@article {Ammari:2020:OSE,
    AUTHOR = {Ammari, Habib and Imeri, Kthim and Nigam, Nilima},
     TITLE = {Optimization of {S}teklov-{N}eumann eigenvalues},
   JOURNAL = {J. Comput. Phys.},
  FJOURNAL = {Journal of Computational Physics},
    VOLUME = {406},
      YEAR = {2020},
     PAGES = {109211, 15},
      ISSN = {0021-9991,1090-2716},
       DOI = {10.1016/j.jcp.2019.109211},
       URL = {https://doi.org/10.1016/j.jcp.2019.109211},
}

@article {Anoop:2018:SMF,
    AUTHOR = {Anoop, T. V. and Bobkov, Vladimir and Sasi, Sarath},
     TITLE = {On the strict monotonicity of the first eigenvalue of the
              {$p$}-{L}aplacian on annuli},
   JOURNAL = {Trans. Amer. Math. Soc.},
  FJOURNAL = {Transactions of the American Mathematical Society},
    VOLUME = {370},
      YEAR = {2018},
    NUMBER = {10},
     PAGES = {7181--7199},
      ISSN = {0002-9947},
       DOI = {10.1090/tran/7241},
       dURL = {https://doi.org/10.1090/tran/7241},
}

@article {Anoop:2021:SVR,
    AUTHOR = {Anoop, T. V. and Ashok Kumar, K. and Kesavan, S.},
     TITLE = {A shape variation result via the geometry of eigenfunctions},
   JOURNAL = {J. Differential Equations},
  FJOURNAL = {Journal of Differential Equations},
    VOLUME = {298},
      YEAR = {2021},
     PAGES = {430--462},
      ISSN = {0022-0396},
       DOI = {10.1016/j.jde.2021.07.001},
       URL = {https://doi.org/10.1016/j.jde.2021.07.001},
}

@article {Arrieta:2008:FRL,
    AUTHOR = {Arrieta, Jos\'{e} M. and Jim\'{e}nez-Casas, \'{A}ngela and
              Rodr\'{\i}guez-Bernal, An\'{\i}bal},
     TITLE = {Flux terms and {R}obin boundary conditions as limit of
              reactions and potentials concentrating at the boundary},
   JOURNAL = {Rev. Mat. Iberoam.},
  FJOURNAL = {Revista Matem\'{a}tica Iberoamericana},
    VOLUME = {24},
      YEAR = {2008},
    NUMBER = {1},
     PAGES = {183--211},
      ISSN = {0213-2230},
       DOI = {10.4171/RMI/533},
}

@book {Bandle:1980:IIA,
    AUTHOR = {Bandle, Catherine},
     TITLE = {Isoperimetric inequalities and applications},
    SERIES = {Monographs and Studies in Mathematics},
    VOLUME = {7},
 PUBLISHER = {Pitman (Advanced Publishing Program), Boston, Mass.-London},
      YEAR = {1980},
     PAGES = {x+228},
      ISBN = {0-273-08423-2},
}

@article {Chorwadwala:2013:TFL,
    AUTHOR = {Chorwadwala, A. M. H. and Vemuri, M. K.},
     TITLE = {Two functionals connected to the {L}aplacian in a class of
              doubly connected domains on rank one symmetric spaces of
              non-compact type},
   JOURNAL = {Geom. Dedicata},
  FJOURNAL = {Geometriae Dedicata},
    VOLUME = {167},
      YEAR = {2013},
     PAGES = {11--21},
      ISSN = {0046-5755},
   MRCLASS = {35J20 (35J05 49R05 58J32)},

       DOI = {10.1007/s10711-012-9800-7},
       dURL = {https://doi.org/10.1007/s10711-012-9800-7},
}

@article {Colbois:2024:SRD,
    AUTHOR = {Colbois, Bruno and Girouard, Alexandre and Gordon, Carolyn and
              Sher, David},
     TITLE = {Some recent developments on the {S}teklov eigenvalue problem},
   JOURNAL = {Rev. Mat. Complut.},
  FJOURNAL = {Revista Matem\'atica Complutense},
    VOLUME = {37},
      YEAR = {2024},
    NUMBER = {1},
     PAGES = {1--161},
      ISSN = {1139-1138,1988-2807},
       DOI = {10.1007/s13163-023-00480-3},
       URL = {https://doi.org/10.1007/s13163-023-00480-3},
}

@book {Colton:1998:IAE,
    AUTHOR = {Colton, David and Kress, Rainer},
     TITLE = {Inverse acoustic and electromagnetic scattering theory},
    SERIES = {Applied Mathematical Sciences},
    VOLUME = {93},
   EDITION = {Fourth},
 PUBLISHER = {Springer, Cham},
      YEAR = {2019},
     PAGES = {xxii+518},
      ISBN = {978-3-030-30350-1; 978-3-030-30351-8},
       DOI = {10.1007/978-3-030-30351-8},
       URL = {https://doi.org/10.1007/978-3-030-30351-8},
}

@incollection {Dittmar:1998:IIS,
    AUTHOR = {Dittmar, Bodo},
     TITLE = {Isoperimetric inequalities for the sums of reciprocal
              eigenvalues},
 BOOKTITLE = {Progress in partial differential equations, {V}ol. 1
              ({P}ont-\`a-{M}ousson, 1997)},
    SERIES = {Pitman Res. Notes Math. Ser.},
    VOLUME = {383},
     PAGES = {78--87},
 PUBLISHER = {Longman, Harlow},
      YEAR = {1998},
      ISBN = {0-582-31708-8},
   MRCLASS = {35P15 (35J25)},

}

@incollection {Dittmar:2005:EPC,
    AUTHOR = {Dittmar, Bodo},
     TITLE = {Eigenvalue problems and conformal mapping},
 BOOKTITLE = {Handbook of complex analysis: geometric function theory.
              {V}ol. 2},
     PAGES = {669--686},
 PUBLISHER = {Elsevier Sci. B. V., Amsterdam},
      YEAR = {2005},
      ISBN = {0-444-51547-X},
       DOI = {10.1016/S1874-5709(05)80020-4},
       URL = {https://doi.org/10.1016/S1874-5709(05)80020-4}
}

@article {Ftouhi:2022:WPS,
    AUTHOR = {Ftouhi, Ilias},
     TITLE = {Where to place a spherical obstacle so as to maximize the
              first nonzero {S}teklov eigenvalue},
   JOURNAL = {ESAIM Control Optim. Calc. Var.},
  FJOURNAL = {ESAIM. Control, Optimisation and Calculus of Variations},
    VOLUME = {28},
      YEAR = {2022},
     PAGES = {Paper No. 6, 21},
      ISSN = {1292-8119},
       DOI = {10.1051/cocv/2021109},
       URL = {https://doi.org/10.1051/cocv/2021109},
}

@article {Gavitone:2022:IIF,
    AUTHOR = {Gavitone, Nunzia and Paoli, Gloria and Piscitelli, Gianpaolo
              and Sannipoli, Rossano},
     TITLE = {An isoperimetric inequality for the first
              {S}teklov-{D}irichlet {L}aplacian eigenvalue of convex sets
              with a spherical hole},
   JOURNAL = {Pacific J. Math.},
  FJOURNAL = {Pacific Journal of Mathematics},
    VOLUME = {320},
      YEAR = {2022},
    NUMBER = {2},
     PAGES = {241--259},
      ISSN = {0030-8730},
   MRCLASS = {35P15 (28A75 35J25)},
 
       DOI = {10.2140/pjm.2022.320.241},
       URL = {https://doi.org/10.2140/pjm.2022.320.241},
}

@book {Gilbarg:2001:EPS,
    AUTHOR = {Gilbarg, David and Trudinger, Neil S.},
     TITLE = {Elliptic partial differential equations of second order},
    SERIES = {Classics in Mathematics},
      NOTE = {Reprint of the 1998 edition},
 PUBLISHER = {Springer-Verlag, Berlin},
      YEAR = {2001},
     PAGES = {xiv+517},
      ISBN = {3-540-41160-7},
}

@article {Girouard:2017:SGS,
    AUTHOR = {Girouard, Alexandre and Polterovich, Iosif},
     TITLE = {Spectral geometry of the {S}teklov problem (survey article)},
   JOURNAL = {J. Spectr. Theory},
  FJOURNAL = {Journal of Spectral Theory},
    VOLUME = {7},
      YEAR = {2017},
    NUMBER = {2},
     PAGES = {321--359},
      ISSN = {1664-039X},
       DOI = {10.4171/JST/164},
       dURL = {https://doi.org/10.4171/JST/164},
}

@article {Hassannezhad:2020:EBM,
    AUTHOR = {Hassannezhad, Asma and Laptev, Ari},
     TITLE = {Eigenvalue bounds of mixed {S}teklov problems},
      NOTE = {With an appendix by Francesco Ferrulli and Jean Lagac\'{e}},
   JOURNAL = {Commun. Contemp. Math.},
  FJOURNAL = {Communications in Contemporary Mathematics},
    VOLUME = {22},
      YEAR = {2020},
    NUMBER = {2},
     PAGES = {1950008, 23},
      ISSN = {0219-1997,1793-6683},
   MRCLASS = {35P15 (35J05 35P20 35S05)},
  
       DOI = {10.1142/S0219199719500081},
       URL = {https://doi.org/10.1142/S0219199719500081},
}

@book {Henrot:2005:VOF,
    AUTHOR = {Henrot, Antoine and Pierre, Michel},
     TITLE = {Variation et optimisation de formes},
    SERIES = {Math\'ematiques \& Applications (Berlin) [Mathematics \&
              Applications]},
    VOLUME = {48},
      NOTE = {Une analyse g\'eom\'etrique. [A geometric analysis]},
 PUBLISHER = {Springer, Berlin},
      YEAR = {2005},
     PAGES = {xii+334},
      ISBN = {978-3-540-26211-4; 3-540-26211-3},
       DOI = {10.1007/3-540-37689-5},
       URL = {https://doi.org/10.1007/3-540-37689-5},
}

@article {Hersch:1968:EPI,
    AUTHOR = {Hersch, Joseph and Payne, Lawrence E.},
     TITLE = {Extremal principles and isoperimetric inequalities for some
              mixed problems of {S}tekloff's type},
   JOURNAL = {Z. Angew. Math. Phys.},
  FJOURNAL = {Zeitschrift f\"{u}r Angewandte Mathematik und Physik. ZAMP.
              Journal of Applied Mathematics and Physics. Journal de
              Math\'{e}matiques et de Physique Appliqu\'{e}es},
    VOLUME = {19},
      YEAR = {1968},
     PAGES = {802--817},
      ISSN = {0044-2275},
   MRCLASS = {35.80},
  
       DOI = {10.1007/BF01591011},
}

@article {Hong:2022:FSD,
    AUTHOR = {Hong, Jiho and Lim, Mikyoung and Seo, Dong-Hwi},
     TITLE = {On the first {S}teklov-{D}irichlet eigenvalue for eccentric
              annuli},
   JOURNAL = {Ann. Mat. Pura Appl. (4)},
  FJOURNAL = {Annali di Matematica Pura ed Applicata. Series IV},
    VOLUME = {201},
      YEAR = {2022},
    NUMBER = {2},
     PAGES = {769--799},
      ISSN = {0373-3114},
       DOI = {10.1007/s10231-021-01137-y},
       URL = {https://doi.org/10.1007/s10231-021-01137-y},
}

@article{Hong:2025:AFS,
    AUTHOR = {Hong, Jiho and Lee, Woojoo and Lim, Mikyoung},
     TITLE = {Approximation of the first {S}teklov-{D}irichlet eigenvalue on
              eccentric spherical shells in general dimensions},
   JOURNAL = {J. Differential Equations},
  FJOURNAL = {Journal of Differential Equations},
    VOLUME = {437},
      YEAR = {2025},
     PAGES = {Paper No. 113295, 24},
      ISSN = {0022-0396,1090-2732},
       DOI = {10.1016/j.jde.2025.113295},
       URL = {https://doi.org/10.1016/j.jde.2025.113295},
}

@incollection {Lamberti:2015:VSE,
    AUTHOR = {Lamberti, Pier Domenico and Provenzano, Luigi},
     TITLE = {Viewing the {S}teklov eigenvalues of the {L}aplace operator as
              critical {N}eumann eigenvalues},
 BOOKTITLE = {Current trends in analysis and its applications},
    SERIES = {Trends Math.},
     PAGES = {171--178},
 PUBLISHER = {Birkh\"{a}user/Springer, Cham},
      YEAR = {2015}
}

@misc{Micetti:2022:SDS,
  title={Steklov-Dirichlet spectrum: stability, optimization and continuity of eigenvalues}, 
      author={Marco Michetti},
      year={2022},
      eprint={2202.08664},
      archivePrefix={arXiv},
      primaryClass={math.AP},
      url={https://arxiv.org/abs/2202.08664},  
}

@article {Paoli:2021:SRS,
    AUTHOR = {Paoli, Gloria and Piscitelli, Gianpaolo and Sannipoli,
              Rossanno},
     TITLE = {A stability result for the {S}teklov {L}aplacian eigenvalue
              problem with a spherical obstacle},
   JOURNAL = {Commun. Pure Appl. Anal.},
  FJOURNAL = {Communications on Pure and Applied Analysis},
    VOLUME = {20},
      YEAR = {2021},
    NUMBER = {1},
     PAGES = {145--158},
      ISSN = {1534-0392},
       DOI = {10.3934/cpaa.2020261},
       URL = {https://doi.org/10.3934/cpaa.2020261},
}

@article {Ramm:1998:IME,
    AUTHOR = {Ramm, Alexander G. and Shivakumar, Pappur N.},
     TITLE = {Inequalities for the minimal eigenvalue of the {L}aplacian in
              an annulus},
   JOURNAL = {Math. Inequal. Appl.},
  FJOURNAL = {Mathematical Inequalities \& Applications},
    VOLUME = {1},
      YEAR = {1998},
    NUMBER = {4},
     PAGES = {559--563},
      ISSN = {1331-4343},
       DOI = {10.7153/mia-01-54},
       dURL = {https://doi.org/10.7153/mia-01-54},
}

@article {Rane:2019:FDE,
    AUTHOR = {Rane, Akanksha V. and Aithal, A. R.},
     TITLE = {The first {D}irichlet eigenvalue of the {L}aplacian in a class
              of doubly connected domains in complex projective space},
   JOURNAL = {Indian J. Pure Appl. Math.},
  FJOURNAL = {Indian Journal of Pure and Applied Mathematics},
    VOLUME = {50},
      YEAR = {2019},
    NUMBER = {1},
     PAGES = {69--81},
      ISSN = {0019-5588},
   MRCLASS = {58J50},
       DOI = {10.1007/s13226-019-0307-x},
      }

@article {Seo:2021:SOP,
    AUTHOR = {Seo, Dong-Hwi},
     TITLE = {A shape optimization problem for the first mixed
              {S}teklov-{D}irichlet eigenvalue},
   JOURNAL = {Ann. Global Anal. Geom.},
  FJOURNAL = {Annals of Global Analysis and Geometry},
    VOLUME = {59},
      YEAR = {2021},
    NUMBER = {3},
     PAGES = {345--365},
      ISSN = {0232-704X},
       DOI = {10.1007/s10455-021-09756-7},
       URL = {https://doi.org/10.1007/s10455-021-09756-7},
}

@article {Stekloff:1902:FPM,
    AUTHOR = {Stekloff, W.},
     TITLE = {Sur les probl\`emes fondamentaux de la physique math\'{e}matique
              (suite et fin)},
   JOURNAL = {Ann. Sci. \'{E}cole Norm. Sup. (3)},
  FJOURNAL = {Annales Scientifiques de l'\'{E}cole Normale Sup\'{e}rieure. Troisi\`eme
              S\'{e}rie},
    VOLUME = {19},
      YEAR = {1902},
     PAGES = {455--490},
      ISSN = {0012-9593}

}

@article {Verma:2020:EPL,
    AUTHOR = {Verma, Sheela and Santhanam, G.},
     TITLE = {On eigenvalue problems related to the laplacian in a class of
              doubly connected domains},
   JOURNAL = {Monatsh. Math.},
  FJOURNAL = {Monatshefte f\"{u}r Mathematik},
    VOLUME = {193},
      YEAR = {2020},
    NUMBER = {4},
     PAGES = {879--899},
      ISSN = {0026-9255},
       DOI = {10.1007/s00605-020-01466-9},
       URL = {https://doi.org/10.1007/s00605-020-01466-9},
}

\end{document}